\documentclass[12pt]{article}
\usepackage{fullpage,epsfig,graphics,amsbsy,amssymb,cancel}
\usepackage{psfrag,hyperref}
\usepackage{graphicx,color}
\usepackage{amsfonts}
\usepackage{amssymb,amsmath,amscd}

\newcommand{\BB}{\mathbb}

\newcommand{\bea}{\begin{eqnarray}}
\newcommand{\eea}{\end{eqnarray}}
\newcommand{\beq}{\begin{eqnarray}}
\newcommand{\eeq}{\end{eqnarray}}
\newcommand{\nn}{\nonumber}

\def\G{\Gamma}
\def\bra{{\langle}}
\def\ket{{\rangle}}
\newcommand{\smalint}{{\Large\textrm{$\smallint$}}}

\newcommand{\FR}{\mathfrak}
\def\ga{\alpha}
\def\gb{\beta}
\def\gc{\gamma}
\def\gd{\delta}

\newtheorem{Thm}{Theorem}[section]

\newtheorem{Def}[Thm]{Definition}
\newtheorem{Exa}[Thm]{Example}

\begin{document}

\thispagestyle{empty}
\begin{flushright} \small
UUITP-14/11\\
 \end{flushright}
\smallskip
\begin{center} \Large
{\bf Introduction to Graded Geometry, \\
Batalin-Vilkovisky Formalism and their Applications}
  \\[12mm] \normalsize
{\bf Jian Qiu$^a$ and Maxim Zabzine$^b$} \\[8mm]
 {\small\it
 ${}^a$I.N.F.N. and Dipartimento di Fisica\\
     Via G. Sansone 1, 50019 Sesto Fiorentino - Firenze, Italy\\
   \vspace{.3cm}
     ${}^b$Department of Physics and Astronomy,
     Uppsala university,\\
     Box 516,
     SE-751\;20 Uppsala,
     Sweden\\
 }
\end{center}
\vspace{7mm}

\begin{abstract}
\noindent
  These notes are intended to provide a self-contained introduction to the basic ideas of  finite dimensional
  Batalin-Vilkovisky (BV) formalism and its applications.  A brief exposition of super- and graded geometries is also
   given.   The BV-formalism is introduced through an odd Fourier  transform and the algebraic aspects of integration
    theory are stressed.
  As a main application we consider the perturbation theory for certain finite dimensional integrals within BV-formalism.
   As an illustration   we present a proof of the isomorphism between the graph complex and the Chevalley-Eilenberg complex of formal
    Hamiltonian vectors fields. We briefly discuss how these ideas can be  extended to the infinite dimensional setting.
    These notes should be accessible to both physicists and mathematicians.
\end{abstract}

\vspace{4cm}
 \begin{center}
 These notes are based on a series of lectures given by the second author
  at the 31th Winter School   ``Geometry and Physics", Czech
Republic, Srni, January 15 - 22, 2011, and are published in Archivum Mathematicum Tomus 47 (2011), 143-199.
\end{center}

\eject
\normalsize
\tableofcontents

\section{Introduction and motivation}

The principal aim of these lecture notes is to present the basic ideas about the  Batalin-Vilkovisky (BV) formalism
 in finite dimensional setting and to elaborate on its application to the perturbative expansion of finite dimensional integrals.
  We try to make these notes self-contained and therefore they include also some background material about
   super and graded geometries, perturbative expansions and graph theory.  We hope that these notes would be
    accessible for math and physics PhD students.

 Originally  the Batalin-Vilkovisky (BV) formalism (named after Igor Batalin and Grigori Vilkovisky, see the original works
  \cite{Batalin:1981jr, Batalin:1984jr}) was introduced in physics as a way of dealing with gauge theories.
    In particular it offers a prescription to perform path integrals of
   gauge theories. In quantum field theory  the path integral is understood as some sort of integral over infinite dimensional
  functional space. Up to now there is no suitable definition of the path integral and in practice all
  heuristic understanding of the path integral is done by mimicking the manipulations of the finite dimensional
  integrals. Thus a proper understanding of  the formal algebraic manipulations with finite  (infinite) dimensional
   integrals is crucial for a better insight to the path integrals.
   Actually nowadays  the  algebraic and combinatorial techniques play a crucial role in dealing with path integral.
    In this context
   the power of BV formalism is that it is able to capture the algebraic properties of the integration and to describe
    the Stokes theorem as some sort of cocyle condition.   The geometrical aspects of BV theory were clarified and
    formalized by Albert Schwarz in \cite{Schwarz:1992nx} and since then it is well-established mathematical subject.

   The idea of these lectures is to present  the algebraic understanding of finite dimensional (super) integrals
   within the framework of BV-formalism and perturbative expansion. Here our intention is to explain the ideas of
    BV formalism in a simplest possible terms and if possible to motivate different formal constructions.
     Therefore instead of presenting many formal definitions and theorems we explain some of the ideas on
      the concrete examples.  At the same time we would like to show the power of BV formalism and thus
       we conclude this note with a highly non-trivial application of BV in finite dimensional setting: the proof
        of the Kontsevich theorem \cite{Kontsevich:Feynamn}  about the  relation between
         graphs and symplectic geometry.

 The outline for the lecture notes is the following. In sections \ref{susy} we briefly review the basic notions from
  supergeometry, in particular $\mathbb{Z}_2$-graded linear algebra, supermanifolds and the integration theory.
   As main examples we discuss the odd tangent and odd cotangent bundles.  In  section \ref{graded} we briefly sketch
    the $\mathbb{Z}$-graded refinement of the supergeometry.  We present a few examples of the graded manifolds.
  In sections \ref{susy} and \ref{graded}  our exposition of super- and graded geometries are quite sketchy.
   We stress the description in terms of local coordinates and avoid many lengthy formal consideration.
   For the full formal exposition of the subject we recommend the recent books \cite{Varadarajan} and \cite{carmeli}.
  In section \ref{oddFT} we introduce the BV structure on the odd cotangent bundle through the odd
   Fourier transformation.  We discuss the integration theory on the odd cotangent bundle and a version of
    the Stokes theorem.
    We stress the algebraic aspects of the integration within BV formalism and explain how the integral gives rise
     to a certain cocycle.  Section \ref{perturbative} provides the basic introduction into the perturbative analysis of
      the finite dimensional integrals. We explain the perturbation theory
        by looking at the specific examples of the integrals in
       $\mathbb{R}^n$ and $\bigoplus\limits_{i=1}^N \mathbb{R}^{2n}$. Also the relevant concepts from the graph theory
        are briefly reviewed and the Kontsevich theorem is stated.  Section \ref{bv-isomorphism} presents the main application of BV formalism to
        the perturbative     expansion of finite dimensional integrals. In particular we present the proof of the
         Kontsevich result \cite{Kontsevich:Feynamn} about
          the isomorphism between  the graph complex and the Chevalley-Eilenberg complex of formal
    Hamiltonian vectors fields. This proof is a simple consequence of the BV formalism and as far as we are aware
     the present form of the proof did not appear anywhere.
          In section \ref{QFT-summary} we outline other application of the present formalism. We briefly discuss the application for   the infinite dimensional setting in the context of quantum field theory.  At the end of the notes there are a few
           Appendices with some technical details and proofs  which we decided not to include in the main text.

\section{Supergeometry}
\label{susy}

The supergeometry extends classical geometry by allowing  odd coordinates,
  which anticommute, in contrast to usual coordinates which commute.
The global objects obtained by gluing such extended coordinate systems, are supermanifolds.
  In this section we briefly review the basic ideas from the supergeometry with the main emphasis on
   the local coordinates. Due to limited time we ignore the sheaf and categorical aspects of supergeometry,
    which are very important for the proper treatment of the subject (see the  books \cite{Varadarajan} and \cite{carmeli}).

\subsection{Idea}

Before going to the formulas and concrete definitions let us say a few general words about the ideas behind
 the super- and graded geometries.  Consider a smooth manifold $M$ and the smooth functions $C^\infty (M)$ over $M$.
 $C^\infty (M)$ is a commutative ring with the point-wise multiplication of the functions and this ring structure contains rich information about
   the original manifold $M$.   The functions which vanish on the fixed region of $M$ form an ideal of this ring and
    moreover the maximal ideals would correspond to the points on $M$. In modern algebraic geometry one replaces
     $C^\infty (M)$ by any commutative ring and the corresponding  ``manifold" $M$ is called scheme. In supergeometry
      (or graded geometry) we replace the commutative ring of functions with supercommutative ring. Thus
supermanifold generalizes the concept of smooth manifold and algebraic schemes to include anticommuting
coordinates. In this sense the super- and graded geometries are conceptually close to the modern algebraic geometry and
 the methods of studying supermanifolds (graded manifold) are variant of those used in the study of schemes.

\subsection{$\mathbb{Z}_2$-graded linear algebra}

The $\mathbb{Z}_2$-graded vector space $V$ over $\mathbb{R}$ (or $\mathbb{C}$) is vector space
 with decomposition
$$V = V_0 \bigoplus V_1~,$$
 where $V_0$ is called even and $V_1$ is called odd.  Any element of $V$ can be decomposed into even
  and odd components. Therefore it is enough to give the definitions for the homogeneous elements. The parity of $v \in V$,
   we denote $|v|$, is defined for the homogeneous element to be $0$ if $v\in V_0$ and $1$ if $v\in V_1$.
   If $\dim V_0=d_0$ and $\dim V_1 = d_1$ then we will adopt the following notation $V^{d_0|d_1}$ and the combination
    $(d_0, d_1)$ is called superdimension of $V$.
    Within the standard use of the terminology $\mathbb{Z}_2$-graded vector space $V$ is the same as superspace.
  All standard constructions from linear algebra (tensor product, direct sum, duality, etc.) carry over to
  $\mathbb{Z}_2$-linear algebra. For example, the morphism between two superspaces is $\mathbb{Z}_2$-grading
   preserving linear map. It is useful to introduce  the parity reversion functor
    which changes the parity of the components of superspace as follows $(\Pi V)_0 = V_1$ and $(\Pi V)_1 = V_0$.
   For example, by $\Pi \mathbb{R}^n$ we mean the purely odd vector space $\mathbb{R}^{0|n}$.

  If $V$ is associative algebra such that the multiplication respects the grading, i.e.
   $|ab|=|a|+|b|$ ($\rm{mod}~ 2$) for homogeneous elements in $V$ then we will call it superalgebra.
   The endomorphsim of superalgebra $V$ is a derivation $D$ of degree $|D|$ if
   \bea
   D(ab)= (Da)b + (-1)^{|D||a|} a (Db)~.\label{derivation-def}
   \eeq
    For any superalgebra we can construct  Lie bracket
     as follows $[a, b] = ab - (-1)^{|a||b|} ba$. By construction this Lie bracket satisfies the following
      properties
     \bea
     [a, b] = -(-1)^{|a||b|} [b,a]~,\label{superLie1}
     \eea
     \bea
     [a,[b,c]] = [[a,b],c] + (-1)^{|a||b|} [b, [a,c]]~.\label{superLie2}
     \eea
     If in general a
      superspace $V$ is equipped with the bilinear bracket $[~,~]$ satisfying the properties (\ref{superLie1}) and (\ref{superLie2})
   then we call it Lie superalgebra. In principle one can define also the odd version of Lie bracket. Namely we can define
  the bracket $[~,~]$ of parity $\epsilon$ such that $|[a,b]|= |a|+|b|+\epsilon$ (${\rm mod} ~2$). This even (odd) Lie bracket
   satisfies the following properties
     \bea
     [a, b] = -(-1)^{(|a|+\epsilon) (|b|+\epsilon)} [b,a]~,\label{oddsuperLie1}
     \eea
     \bea
    [a, [b, c]]= [[a, b], c] + (-1)^{(|a|+\epsilon)(|b|+\epsilon)} [b, [a,c]]~ .\label{oddsuperLie2}
     \eea
 However the odd  Lie superbracket can be mapped to even Lie superbracket by the parity reversion functor.
 Thus odd case can be always reduced to the even.

  Coming back to the general superalgebras.
  The supergalgebra $V$ is called supercommutative if
     $$ ab = (-1)^{|a||b|} ba~.$$
      The supercommutative algebras will play the central role in our considerations. Let us discuss a very important
       example of the suprecommutative algebra, the exterior algebra.
\begin{Exa}\label{exterioralgebra}
 Consider purely odd superspace   $\Pi \mathbb{R}^m = \mathbb{R}^{0|m}$ over the real number of dimension $m$.
 Let us pick up the basis $\theta^i$, $i=1,2,..., m$ and define the multiplication between the basis elements
  satisfying  $\theta^i \theta^j = - \theta^j \theta^i$. The functions $C^\infty (\mathbb{R}^{0|m})$
  on $\mathbb{R}^{0|m}$  are given by the following expression
 $$f(\theta^1, \theta^2, ..., \theta^m) = \sum\limits_{l=0}^m\frac{1}{l!} ~f_{i_1 i_2 ... i_l}  \theta^{i_1} \theta^{i_2} ... \theta^{i_l}~,$$
  and they correspond to the elements of exterior algebra $\wedge^\bullet (\mathbb{R}^m)^*$.  The exterior algebra
  $$\wedge^\bullet (\mathbb{R}^m)^* = \left (\wedge^{\rm even} ( \mathbb{R}^m)^* \right ) \bigoplus \left (  \wedge^{\rm odd} (\mathbb{R}^m)^*
   \right )$$
    is a supervector space with the supercommutative multiplications given by wedge product. The wedge product of the exterior
     algebra corresponds to the function multiplication in  $C^\infty (\mathbb{R}^{0|m})$.
 \end{Exa}

  Let us consider the supercommutative algebra $V$ with the multiplication and in addition there is
   a Lie bracket of parity $\epsilon$.  We require that  ${\rm ad}_a = [a,~]$ is a derivation of $\cdot$
    of degree $|a|+\epsilon$, namely
     \bea
     [ a, bc] = [a, b]c + (-1)^{(|a|+\epsilon)|b|} b [a,c]~. \label{Poisson-multiplication}
  \eea
   Such structure $(V, \cdot, [~,~])$ is called even Poisson algebra for $\epsilon=0$ and Gerstenhaber algebra
    (odd Poisson algebra) for $\epsilon=1$. It is crucial that it is not possible to reduce Gerstenhaber algebra
     to even Poisson algebra by the parity reversion, since now we have two operations in the game, supercommutative product and
      Lie bracket compatible in a specific way.

\subsection{Supermanifolds}

 We can construct more complicated examples of the supercommutative algebras.
  Consider the real superspace $\mathbb{R}^{n|m}$ and we define the space of functions on it as follows
  $$ C^\infty (\mathbb{R}^{n|m}) \equiv C^\infty (\mathbb{R}^n) \otimes \wedge^\bullet (\mathbb{R}^m)^*~.$$
   If we pick up an open subset $U_0$ in $\mathbb{R}^n$ then we can associate to $U_0$ the supercommutative
 algebras as follows
 \bea
  U_0~\longrightarrow~ C^\infty (U_0) \otimes \wedge^\bullet (\mathbb{R}^m)^*~.\label{freegenerated}
 \eea
  This supercommutative algebra can be thought of as the algebra of functions on the superdomain  $U^{n|m} \subset \mathbb{R}^{n|m}$,
 $C^\infty (U^{n|m}) = C^\infty (U_0) \otimes \wedge^\bullet (\mathbb{R}^m)^*$.
 The   superdomain $U^{n|m} \subset \mathbb{R}^{n|m}$ can be characterized in terms of standard even
    coordinates $x^\mu$ ($\mu=1,2,..., n$) for $U_0$ and the odd coordinates $\theta^i$ ($i=1,2,...,m$), such that $\theta^i \theta^j
     = -\theta^j \theta^i$.  In analogy with ordinary manifolds a supermanifold can be defined by gluing together superdomains
      by degree preserving maps. Thus the domain $U^{n|m}$ with coordinates $(x^\mu, \theta^i)$  can be glued to the domain
       $V^{n|m}$ with coordinates $(\tilde{x}^\mu, \tilde{\theta}^i)$ by invertible and degree-preserving maps  $\tilde{x}^\mu =
        \tilde{x}^\mu (x, \theta)$   and $\tilde{\theta}^i =  \tilde{\theta}(x, \theta)$ defined for $x \in U_0 \cap V_0$. Thus formally  the theory
         of supermanifolds mimics the standard smooth manifolds. However one should anticipate that some of the geometric intuition
          fails and we cannot think in terms of points due to the presence of the odd coordinates.  This situation is very similar
           to the algebraic geometry when there can be nilpotent elements in the commutative ring.

   The supermanifold is defined by gluing superdomains. However, the gluing should be done with some care and
  for the rigorous treatment we need to use the sheaf theory. Let us give a precise definition of the smooth supermanifold.
  \begin{Def}
  A smooth supermanifold ${\cal M}$ of dimension $(n,m)$
   is a smooth manifold $M$ with a sheaf of supercommutative superalgebras, typically denoted $O_{M}$ or $C^\infty ({\cal M})$, that
    is locally isomorphic to $C^\infty (U_0) \otimes \wedge^\bullet (\mathbb{R}^m)^*$, where $U_0$ is open subset of $\mathbb{R}^n$.
  \end{Def}
  Thus essentially the supermanifold is defined through the gluing supercommutative algebras which locally look like in
   (\ref{freegenerated}).  This supercommutative algebra is sometimes called 'freely generated' since it can be generated by
    even and odd coordinates $x^\mu$ and $\theta^i$. If we allow more general supercommutative algebras to be glued, we
     will be led to the notion of superscheme which is a natural super generalization in the algebraic geometry.

Let us illustrate this formal definition of supermanifold with couple of concrete examples.

\begin{Exa}\label{oddtangentbudle} Assume that $M$ is smooth manifold then we can associate to it the supermanifold
$\Pi TM$ odd tangent bundle, which is defined by the gluing rule
$$\tilde{x}^\mu =\tilde{x}^\mu (x)~,~~~~~~~\tilde{\theta}^\mu = \frac{\partial \tilde{x}^\mu}{\partial x^\nu} \theta^\nu~,$$
  where $x$'s are local coordinates on $M$ and $\theta$'s are glued as $dx^\mu$. The functions on $\Pi TM$ have the following
   expansion
$$ f (x, \theta) = \sum\limits_{p=0}^{\dim M} \frac{1}{p!} f_{\mu_1\mu_2 ... \mu_p} (x) \theta^{\mu_1} \theta^{\mu_2} ... \theta^{\mu_p}$$
 and thus they are naturally identified with the differential forms,  $C^\infty (\Pi TM)= \Omega^\bullet (M)$. Indeed
  locally the differential forms correspond to freely generated supercommutative algebra
  $$ \Omega^\bullet (U_0) = C^\infty (U_0)  \otimes \wedge (\mathbb{R}^n)^*~.$$
 \end{Exa}

\begin{Exa}\label{oddcotangentbundle} Again let $M$ be a smooth manifold and we associate to it now another super manifold
$\Pi T^*M$ odd cotangent bundle, which has the following local description
$$\tilde{x}^\mu =\tilde{x}^\mu (x)~,~~~~~~~\tilde{\theta}_\mu =  \frac{\partial x^\nu}{\partial \tilde{x}^\mu} \theta_\nu ~,$$
 where $x$'s are local coordinates on $M$ and $\theta$'s transform as $\partial_\mu$. The functions on $\Pi T^*M$ have the expansion
$$ f (x, \theta) = \sum\limits_{p=0}^{\dim M} \frac{1}{p!} f^{\mu_1\mu_2 ... \mu_p} (x) \theta_{\mu_1} \theta_{\mu_2} ... \theta_{\mu_p}$$
 and thus they are naturally identified with multivector fields,  $C^\infty (\Pi T^*M) = \Gamma (\wedge^\bullet TM)$. Indeed
  the sheaf of multivector fields is a sheaf of supercommutative algebras which is locally freely generated.
\end{Exa}

  Many notions and results from the standard differential geometry can be extended to supermanifolds in straightforward fashion.
  For example, the vector fields on supermanifold ${\cal M}$ are defined as derivations of the supercommutative algebra $C^\infty ({\cal M})$.
  The use of local coordinates is extremely powerful and sufficient for most purposes. The notion of morphisms of supermanifolds can be described
   locally exactly as it is done in the case of smooth manifolds.

\subsection{Integration theory}
\label{subs-integration}

Now we have to discuss the integration theory for the supermanifolds. We need to define the measure and it can be done
 first locally in analogy with the standard case.  The main novelty comes from the odd part of the measure.

Let us start from the discussion of the integration of the function  $f(x)$ in one variable. The even
integral is defined as usual
\bea
 \int f(x) dx
\eea
 and if we change the coordinate $\tilde{x} = c x$ then the measure is changed accordingly to the standard rules $d\tilde{x} = c dx$.
 Next consider the function of one odd variable $\theta$ which is given by $f = f_0 + f_1 \theta$, where $f_0$ and $f_1$ are some
  real numbers.  We define the integral over this function as linear operation such that
 \bea
  \int d \theta =0~,~~~~~~~~~~\int d\theta~ \theta = 1~.
 \eea
    Now  if we change the odd coordinate $\tilde{\theta} = c \theta$ we still want the same
     definition to hold, namely
   \bea
  \int d \tilde{\theta} =0~,~~~~~~~~~~\int  d\tilde{\theta}~ \tilde{\theta}= 1~.
 \eea
  As a result of this  we get  that  the odd measure transforms as
   follows $d\tilde{\theta} = \frac{1}{c} d\theta$ and this transformation property should be contrasted with the even integration.
   Next we can define the odd measure over functions of many $\theta$'s.  Assume that  there are odd $\theta^i$ ($i=1,2,..., m$).
    Using the definition for a single $\theta$ we define the measure to be such that
   \bea
\int d^m \theta~ \theta^1\theta^2 ... \theta^m \equiv    \int d\theta^n ... \int  d\theta^2 \int d\theta^1 ~\theta^1 \theta^2 ... \theta^m
  = 1
   \eea
    and all other integrals are zero.  Let us change variables according to the following rule
   $\tilde{\theta}^i = A^i_j \theta^j$ such that
   $$ \tilde{\theta}^1 \tilde{\theta}^2~ ...~ \tilde{\theta}^m = \det A~ \theta^1 \theta^2 ~... ~\theta^m~.$$
    In new variables  we still require that
  \bea
    \int d^n \tilde{\theta}~ \tilde{\theta}^1\tilde{\theta}^2 ... \tilde{\theta}^n =1~.
\eea
 Therefore  we obtain the following formula for the transformation of the measure,  $d^n \tilde{\theta} = (\det A)^{-1} d^n \theta$.
  Using these simple ideas we can define the integration of the function over any superdomain $U^{n|m}$ and then we have
   to check how the measure is glued as we patch different superdomains.
  On a supermanifold we would like to integrate the functions and for this we will need well-defined measure of
  the integration on the whole supermanifold.

  Instead of writing down the general formulas let us discuss the integration of functions on odd tangent and odd cotangent bundles.

\begin{Exa}\label{integrationtangent}
 Using the notation from the example \ref{oddtangentbudle} let us study the integration measure on
  the odd tangent bundle $\Pi TM$. The even part of the measure transforms in the standard way
 $$ d^n \tilde{x} = \det \left ( \frac{\partial \tilde{x}}{\partial x} \right ) d^n x~,$$
 while the odd part transforms according to the following property
 $$ d^n \tilde{\theta} =  \frac{1}{\det \left ( \frac{\partial \tilde{x}}{\partial x} \right )} d^n \theta~.$$
  As we can see the transformation of even and odd parts cancel each other and thus we have
  $$ \int d^n\tilde{x}~ d^n \tilde{\theta} = \int d^n x~ d^n \theta~,$$
 which corresponds to the canonical integration on $\Pi TM$. Any function of top degree on $\Pi TM$ can be integrated
  canonically. This is not surprising since the integration of the top differential forms is defined canonically for any smooth orientable
   manifold.
 \end{Exa}
\begin{Exa}\label{integrationcotangent} Using the notation from the example \ref{oddcotangentbundle} let us study
the integration on odd cotangent bundle $\Pi T^*M$. The even part transforms as before
$$ d^n \tilde{x} = \det \left ( \frac{\partial \tilde{x}}{\partial x} \right ) d^n x~,$$
 while the odd part transforms in the same way as even
 $$ d^n \tilde{\theta} = \det  \left ( \frac{\partial \tilde{x}}{\partial x} \right )d^n \theta~.$$
Assume that $M$ is orientiable and  let us pick up a volume form (nowhere vanishing top form)
 $${\rm vol} = \rho (x)~ dx^1 \wedge ... \wedge dx^n~.$$
  One can check that $\rho$ transforms as a densitity
   $$\tilde{\rho} =  \frac{1}{\det \left ( \frac{\partial \tilde{x}}{\partial x} \right )} \rho~. $$
    Combining all these ingredients together we can define the following invariant measure
     $$\int d^n \tilde{x}~ d^n \tilde{\theta} ~\tilde{\rho}^2= \int d^n x ~d^n \theta ~ \rho^2~,$$
      which we can glue consistently.
   Thus to integrate the multivector fields we need to pick a volume form on $M$.
 \end{Exa}

\section{Graded geometry}
\label{graded}

Graded geometry is $\mathbb{Z}$-refinement of supergeometry. Many definitions from the supergeometry have
 straightforward generalization to the graded case.
  In our review of graded geometry we will be very brief,  for more details one can consult \cite{Roytenberg:2002nu,susy-note-BB}.

\subsection{$\mathbb{Z}$-graded linear algebra}

A $\mathbb{Z}$-graded vector space is a vector space $V$ with the decomposition labelled by integers
$$V = \bigoplus\limits_{i \in \mathbb{Z}} V_i~.$$
 If  $v \in V_i$ then we say that $v$ is homogeneous element of $V$ a degree $|v|=i$. Any element of $V$ can
  be decomposed in terms of homogeneous elements of a given degree.  Many concepts of linear algebra and superalgebra
    has a straightforward generalization to the general graded case.  The morphism between graded vector spaces is defined as
    a linear map which preserves the grading.  Assuming that $\mathbb{R}$ (or $\mathbb{C}$) is vector space of degree $0$
    the dual vector space $(V_i)^*$  is defined as  $V_{-i}^*$.
    The graded vector space $V[k]$ shifted by degree $k$ is defined as direct sum of $V_{i+k}$.

  If the graded vector space $V$ is equipped with the associative product which respects the grading then we call
   $V$  a graded algebra.  The endomorphism of graded algebra $V$ is a derivation $D$ of degree $|D|$ if it satisfies the relation (\ref{derivation-def}), but now with $\mathbb{Z}$-grading.
   If for a graded algebra $V$ and any homogeneous elements $v$ and $\tilde{v}$ therein we have the relation
    $$ v \tilde{v}  = (-1)^{|v| |\tilde{v}|} \tilde{v}  v~,$$
then we call $V$ a graded commutative algebra. The graded commutative algebras play the crucial role in the graded geometry.
 One of the most important examples of graded algebra is given by the graded symmetric space $S(V)$.

  \begin{Def}\label{def-symmalg-graded}
   Let $V$ be a graded vector space over $\BB{R}$ or $\BB{C}$. We define the
graded symmetric algebra $S(V)$  as the linear space spanned by polynomial functions on $V$
$$ \sum\limits_l f_{a_1 a_2 ... a_l}~ v^{a_1} v^{a_2} ... v^{a_l}~,$$
where we use the relations
$$ v^{a} v^{b}  = (-1)^{|v^a| |v^b|}  v^b  v^a$$
 with $v^a$ and $v^b$ being
  homogeneous elements of degree $|v^a|$ and $|v^b|$ respectively.
   The functions on $V$ are naturally graded and multiplication of functions
    is graded commutative.
Therefore the graded symmetric algebra $S(V)$ is a graded commutative algebra.
  \end{Def}

In analogy with $\mathbb{Z}_2$-case we can define the Lie bracket $[~,~]$ of the integer degree $\epsilon$ now such that
 $|[v,w]| = |v| + |w| +\epsilon$ and it satisfies the properties (\ref{oddsuperLie1}) and (\ref{oddsuperLie2}).
   Analogously we can introduce the graded versions
  of Poisson algebra. If the $\mathbb{Z}$-graded vector space $V$ is equipped with a graded commutative algebra structure $\cdot$
   and a Lie algebra bracket $[~,~]$  of degree $\epsilon$ such that they are compatible with respect to the relation (\ref{Poisson-multiplication})
    then we call $V$ $\epsilon$-graded Poisson algebra (or simply $\epsilon$-Poisson algebra).  The standard use of terminology
     is the following, $0$-graded Poisson algebra is  called Poisson algebra and $(\pm 1)$-graded Poisson algebra is called quite
      often Gerstenhaber algebra.  For more explanation and examples of graded Poisson algebras the reader may
       consult \cite{cattaneo-GP}.

 Let us make one important side remark about the sign conventions in the graded case. Quite often one has to deal with
  bi-graded vector spaces which carry simultaneously $\mathbb{Z}_2$- and $\mathbb{Z}$-gradings. There exist two different
   sign conventions when one moves one element past another,
   \bea v w  = (-1)^{pq+ls} w  v~,\label{convention_1}\eea
   and
\bea   v w  = (-1)^{(p+q)(l+s)} w  v~,\label{convention_2}\eea
 where the degrees are defined as follows
  $$ |v|_{\mathbb{Z}_2} =p~,~~~ |v|_{\mathbb{Z}} =l~,~~~|w|_{\mathbb{Z}_2} =q~,~~~|w|_{\mathbb{Z}} =s~.$$
   Both conventions are widely used and they each have their advantages. They are equivalent, but one should never mix them
    while dealing the $\mathbb{Z}$-graded superspaces. For more details see the explanation in  \cite{susy-note-BB}.
     However this sign subtlety is irrelevant for most of our consideration.

\subsection{Graded manifold}

 We can define the graded manifolds very much in analogy with the supermanifolds. We have sets of the coordinates
 with assignment of degree and we glue them by the degree preserving maps.  Let us give the formal definition first.

  \begin{Def}
  A smooth graded manifold ${\cal M}$
   is a smooth manifold $M$ with a sheaf of graded commutative algebras, typically  denoted by $C^\infty ({\cal M})$, which
    is locally isomorphic to $C^\infty (U_0) \otimes S(V)$, where $U_0$ is open subset of $\mathbb{R}^n$ and $V$ is graded
     vector space.
  \end{Def}

  This definition is a generalization of supermanifold to the graded case. To every patch we associate a commutative
   graded algebra which is freely generated by the graded coordinates. The gluing is done by the degree preserving maps.
    The best way of explaining this definition is by considering the explicit examples.

\begin{Exa}\label{gradedtangentbudle}
  Let us introduce the graded version of the odd tangent bundle from the example \ref{oddtangentbudle}. We denote
   the graded tangent bundle as  $T[1]M$ and we have the same coordinates $x^\mu$ and $\theta^\mu$ as in
    the example \ref{oddtangentbudle}, with the same transformation rules. The coordinate $x$ is of degree $0$ and $\theta$ is
     of degree $1$ and the gluing rules respect the degree. The space of functions $C^\infty (T[1]M)=\Omega^\bullet(M)$ is a graded
      commutative  algebra
      with the same $\mathbb{Z}$-grading as the differential forms.
 \end{Exa}

 \begin{Exa}\label{gradedcotangentbudle}
Analogously we can introduce the graded version $T^*[-1]M$ of the odd cotangent bundle from the example \ref{oddcotangentbundle}.
 Now we allocate the degree $0$ for $x$ and degree $-1$ for $\theta$.  The gluing preserves the degrees.
 The functions
  $C^\infty (T^*[-1]M) = \Gamma (\wedge^\bullet TM)$ is graded commutative algebra with degree given by minus of
   degree of multivector field.
 \end{Exa}

 \begin{Exa}\label{gradedcocotangentbudle}
  Let us discuss a slightly more complicated example of graded cotangent bundle over cotangent bundle
  $T^*[2]  (T^*[1] M)$. In local coordinates we can describe it as follows. Introduce the coordinates
   $x^\mu$, $\theta^\mu$, $\psi_\mu$ and $p_\mu$ of degree $0$, $1$, $1$ and $2$ respectively.
    The gluing between patches is done by the following degree preserving maps
  $$\tilde{x}^\mu = \tilde{x}^\mu (x)~,~~~~~\tilde{\theta}^\mu = \frac{\partial \tilde{x}^\mu}{\partial x^\nu} \theta^\nu~,~~~~~
  \tilde{\psi}_\mu = \frac{\partial x^\nu}{\partial \tilde{x}^\mu} \psi_\nu~,$$
  $$ \tilde{p}_\mu = \frac{\partial x^\nu}{\partial \tilde{x}^\mu} p_\nu + \left (\frac{\partial^2 x^\nu}{\partial \tilde{x}^\gamma
   \tilde{x}^\mu} \right ) \frac{\partial \tilde{x}^\gamma}{\partial x^\sigma} \psi_\nu \theta^\sigma~.$$
    Now it is bit more complicated to describe the functions $C^\infty (T^*[2]  (T^*[1] M))$ in terms of standard geometrical
     objects. However by construction $C^\infty (T^*[2]  (T^*[1] M))$ is a graded commutative algebra. In degree zero
      $C^\infty (T^*[2]  (T^*[1] M))$ corresponds to $C^\infty(M)$ and in degree one to $\Gamma (TM \oplus T^*M)$.
       For more details of this example the reader may consult \cite{Roytenberg:2002nu}.
  \end{Exa}

 Again the big  chunk of differential geometry has a straightforward generalization to the graded manifolds.
  The integration theory for the graded manifolds is totally analogous to the super case, with the main
   difference between the even and odd measure described in subsection \ref{subs-integration}.
 The vector fields are defined as derivations of $C^\infty({\cal M})$ for the graded manifold ${\cal M}$.  The vector
  fields on ${\cal M}$ are naturally graded, and amongst these we are interested in the odd vector fields which square to zero.

    \begin{Def}
  If the graded manifold ${\cal M}$ is equipped with a derivation $D$ of $C^\infty ({\cal M})$ of
   degree $1$ with additional property $D^2=0$ then we call such $D$
  a homological vector field.
   $D$ endows the graded commutative algebra of function  $C^\infty ({\cal M})$ with
    the structure of differential complex.  One calls such graded commutative algebra with $D$ a
     graded differential algebra.
  \end{Def}

  Let us state the most important example of homological vector field for the graded tangent bundle.

 \begin{Exa}\label{deRhamdifferential}
  Consider the graded tangent bundle $T[1]M$ described in the example \ref{gradedtangentbudle}.
   Let us introduce the vector field of degree $1$ written in local coordinates as follows
 $$D= \theta^\mu \frac{\partial}{\partial x^\mu}~,$$
  which is glued in an obvious way. Since $D^2=0$ this is an example of homological vector field.
   $D$ on $C^\infty(T[1]M)= \Omega^\bullet (M)$ corresponds to the de Rham differential on $\Omega^\bullet (M)$.
 \end{Exa}

\section{Odd Fourier transform and BV-formalism}
\label{oddFT}

In this section we introduce the basics of BV formalism. We derive the construction
 through the odd Fourier transformation which maps $C^\infty (T[1]M)$ to $C^\infty(T^*[-1]M)$.
   Odd cotangent bundle $T^*[-1]M$ has a nice algebraic structure on the space of functions
    and using the odd Fourier transform we will derive  the version of Stokes theorem for  the
     integration on $T^*[-1]M$.  The power of  BV formalism  is based on the algebraic interpretation
      of the integration theory for odd cotangent bundle.

\subsection{Standard Fourier transform}

Let us start by recalling the well-known properties of the standard Fourier transformations.
Consider the suitable function $f(x)$ on the real line $\mathbb{R}$ and define the
 Fourier transformation of  this function according to the following formula
\bea\label{normalFT1}
F[f] (p) = \frac{1}{\sqrt{2\pi}} \int\limits_{-\infty}^{\infty} f(x) e^{-ipx} dx~.
\eea
 One can also define the inverse Fourier transformation as follows
\bea\label{normalFT2}
 F^{-1}[f] = \frac{1}{\sqrt{2\pi}} \int\limits_{-\infty}^{\infty} f(p) e^{ipx} dp~.
\eea
 There are some subtleties related to the proper understanding of the integrals (\ref{normalFT1})-(\ref{normalFT2})
  and certain restrictions on $f$ to make sense of these expressions. However, let us put aside these complications in this note. The functions
   have associative point-wise multiplication and one can study how it is mapped under the Fourier transformation.
    It is an easy exercise to show that
\bea
 F[f] F[g] = F[ f* g]
\eea
 where $*$-product is defined as follows
 \bea
  (f *g) (x) = \int\limits_{-\infty}^{\infty} f(y) g(x-y) dy~.
 \eea
  This $*$-operation  is called convolution of two functions and it can be defined for any two integrable
   functions on the line.
 This $*$-product is associative $(f*g)*h = f*(g*h)$ and commutative $f*g=g*f$.  Thus the space of integrable functions
  is associative commutative algebra with respect to convolution, but there is no identity (since $1$ is not an integrable function on the
   line). It is important to stress that the derivative $\frac{d}{dx}$ is not a derivation of this $*$-product.

\subsection{Odd Fourier transform}
\label{ss-odd-FT}

 Let us assume that the manifold $M$ is orientable and we can pick up a volume form
\bea
{\rm vol} = \rho (x) ~dx^1 \wedge... \wedge dx^n= \frac{1}{n!}~ \Omega_{\mu_1 ... \mu_n}(x)~dx^{\mu_1} \wedge...
 \wedge dx^{\mu_n}~,\label{defin-volume}
 \eea
 which is a top degree nowhere vanishing form and $n=\dim M$.
  Consider the graded manifold $T[1]M$ and the integration theory
 which we have discussed in  the example \ref{integrationtangent}. If we have the volume form then
  we can define the integration only along the odd direction as follows
  $$ \int d^n \tilde{\theta}~ \tilde{\rho}^{-1} = \int d^n \theta ~\rho^{-1}~.$$
 In analogy with the standard Fourier transform (\ref{normalFT1}) we can define the odd Fourier transfrom
for $f(x, \theta) \in C^\infty (T[1]M)$ as
\bea
 F[f](x, \psi) = \int d^n\theta~\rho^{-1} e^{\psi_\mu \theta^\mu}  f(x, \theta)~,\label{F-T}
\eea
where $d^d{\theta}=d\theta^d\cdots d\theta^1$.
  Obviously  we would like to make sense globally of the transformation (\ref{F-T}). Therefore we assume that
   the degree of $\psi_\mu$ is $-1$ and it transforms as $\partial_\mu$ (so in the way dual to $\theta^\mu$).
   Thus $ F[f](x, \psi) \in C^{\infty} (T^*[-1]M)$ and the odd Fourier transform maps functions on $T[1]M$ to
 the functions on $T^*[-1]M$. The explicit formula (\ref{FT-explicit}) of the Fourier transform of  $p$-form is given
 in the Appendix. We can also define the inverse Fourier transform $F^{-1}$ which maps the functions on $T^*[-1]M$
  to the functions  on $T[1]M$
  as follows
 \bea
  F^{-1} [\tilde{f}] (x, \theta) = (-1)^{n(n+1)/2}\int d^n\psi~\rho~ e^{-\psi_{\mu}\theta^{\mu}}\tilde f(x,\psi)~,\label{inverse-FT}
 \eea
  where $\tilde{f}(x, \psi)  \in C^\infty (T^*[-1]M)$.
 One may easily check  that
 \bea (F^{-1}F[f])(x,\eta)&=&(-1)^{n(n+1)/2}\int d^n\psi~\rho e^{-\psi_{\mu}\eta^{\mu}}~\int d^n\theta~\rho^{-1}~ e^{\psi_{\mu}\theta^{\mu}} f(x,\theta) =f(x,\eta)~.\nn
 \eea
 Since we have to discuss both odd tangent and odd cotangent bundles simultaneously, in this section we adopt the following notation for the functions: we denote with symbols without tilde functions on $T[1]M$ and with tilde functions on $T^*[-1]M$.

$C^\infty(T[1]M)$ is a differential graded algebra with the graded commutative multiplication and the differential $D$ defined in
 example \ref{deRhamdifferential}. Let us discuss how these operations behave  under the odd Fourier transform $F$.
 Under $F$   the differential $D$  transforms to bilinear operation $\Delta$ on $C^\infty(T^*[-1]M)$ as follows
 \bea
  F[Df]=  (-1)^n \Delta F[f] \label{D-Delta-relation}
 \eea
  and from this we can calculate the explicit form of $\Delta$
 \bea
  \Delta = \rho^{-1} \frac{\partial^2}{\partial x^\mu \partial \psi_\mu} \rho = \frac{\partial^2}{\partial x^\mu \partial \psi_\mu}  + \partial_\mu (\log \rho)
  \frac{\partial}{\partial \psi_\mu}~ .\label{def-Laplace}
 \eea
  By construction $\Delta^2=0$ and degree of $\Delta$ is $1$.  Next let us discuss how the graded commutative product on
   $C^\infty(T[1]M)$ transforms under $F$.  The situation is very much analogous to the standard Fourier transform where
    the multiplication of functions goes to their convolution.
     To be specific we have
 \bea
F[f g] = F[f]*F[g]
\eea
and from this we derive the explicit formula for the odd convolution
\bea
 (\tilde{f}*\tilde{g})(x, \psi) = (-1)^{n(n+|f|)} \int d^n\lambda~ \rho~  \tilde{f}(x, \lambda) \tilde{g}(x, \psi- \lambda)~,
\eea
 where $\tilde{f}, \tilde{g} \in C^\infty (T^*[-1]M)$ and $\psi$, $\lambda$ are odd coordinates on $T^*[-1]M$.
  This star product is associative and by construction $\Delta$ is a
  derivation of this product (since $D$ is a derivation of usual product on $C^\infty(T[1]M)$). Moreover we have
   the following relation
 \bea
  \tilde{f} * \tilde{g} = (-1)^{(n-|\tilde{f}|)(n-|\tilde{g}|)} \tilde{g}* \tilde{f}
 \eea
  and thus this star product does not preserve $\mathbb{Z}$-grading, i.e. $|\tilde{f}* \tilde{g}| \neq |\tilde{f}| + |\tilde{g}|$.
 Thus the odd convolution  of functions is not a graded commutative product, which should not be surprising since
  $F$ is not a morphism of the graded manifolds (generically it is not a morphisms of supermanifolds either).
   At the same time $C^\infty(T^*[-1]M)$ is a graded commutative algebra with respect to the ordinary multiplication
    of functions, but $\Delta$ is not a derivation of this multiplication
 \bea
  \Delta(\tilde{f}\tilde{g}) \neq \Delta(\tilde{f}) \tilde{g} +  (-1)^{|\tilde{f}|} \tilde{f} \Delta (\tilde{g})~.
 \eea
 We can define the bilinear operation which measures the failure of $\Delta$ to be a derivation
  \bea
(-1)^{|\tilde{f}|} \{\tilde{f}, \tilde{g}\}=  \Delta(\tilde{f}\tilde{g}) - \Delta(\tilde{f}) \tilde{g} - (-1)^{|\tilde{f}|} \tilde{f} \Delta (\tilde{g})~.\label{Delta-defPB}
 \eea
A direct calculation gives the following expression
 \bea
  \{ \tilde{f}, \tilde{g}\} = \frac{\partial \tilde{f}}{\partial x^\mu} \frac{\partial \tilde{g}}{\partial \psi_\mu} + (-1)^{|f|} \frac{\partial \tilde{f}}{\partial \psi_\mu} \frac{\partial \tilde{g}}{\partial x^\mu}~,\label{def-PBoncot}
 \eea
  which is very reminiscent of the standard Poisson bracket for the cotangent bundle, but now with the odd momenta.
   For the derivative $\frac{\partial }{\partial \psi_\mu}$ we use the following convention
   \bea \frac{\partial \psi_\nu}{\partial \psi_\mu} = \delta_\nu^\mu \nn \eea
   and it is derivation of degree $1$ (see the definition (\ref{derivation-def})).
  By a direct calculation
   one can check that this bracket (\ref{def-PBoncot}) gives rise to $1$-Poisson algebra (Gerstenhaber algebra) on $C^\infty(T^*[-1]M)$.
    Indeed the bracket (\ref{def-PBoncot}) on $C^\infty(T^*[-1]M)$ corresponds to the Schouten bracket on the multivector fields (see Appendix
    for the explicit formulas).
    To summarize, upon the choice of  volume form on $M$, $C^\infty(T^*[-1]M)$ is an odd Poisson algebra (Gerstenhaber algebra) with
     the Poisson bracket generated by $\Delta$-operator as in (\ref{Delta-defPB}). Such a structure is called BV-algebra. We
      will now summarize and formalize this notion.

Let us recall the definition of odd Poisson algebra (Gerstenhaber algebra).
\begin{Def}\label{Gerstenhaber-algebra}
The graded commutative algebra $V$ with the odd bracket $\{~,~\}$ satisfying the following axioms
$$ \{v, w\} = - (-1)^{(|v|+1)(|w|+1)} \{w,v\}$$
$$ \{v, \{w, z\}\} = \{\{v, w\}, z\} + (-1)^{(|v|+1)(|w|+1)} \{w, \{v, z \}\}$$
 $$ \{ v, wz\} = \{v, w\} z + (-1)^{(|v|+1)|w|} w \{v, z\}$$
  is called a Gerstenhaber algebra.
\end{Def}
 Typically it is assumed that the degree of  bracket $\{~,~\}$ is 1 (or $-1$ depending on conventions).
 Thus the space of functions $C^{\infty} (T^*[-1]M)$ is a Gerstenhaber algebra with a graded commutative multiplication of
  functions and a bracket of degree $1$ defined by (\ref{def-PBoncot}).
  The BV-algebra is Gerstenhaber algebra with an additional structure.
\begin{Def}\label{BValgebra}
 A Gerstenhaber algebra $(V, \cdot,  \{~,~\})$ together with an odd $\mathbb{R}$--linear map
 $$\Delta ~:~V \longrightarrow V~,$$
  which squares to zero $\Delta^2=0$ and generates the bracket $\{~,~\}$ according to
  \bea
   \{v, w\}=  (-1)^{|v|} \Delta(vw)+(-1)^{|v|+1}(\Delta v)w - v(\Delta w)~,\label{def-bracket-Delta}
   \eea
  is called a BV-algebra. $\Delta$ is called the odd Laplace operator (odd Laplacian).
 \end{Def}
Again it is assumed that degree $\Delta$ is 1 (or $-1$ depending on conventions).
The space of functions $C^\infty (T^*[-1]M)$ is a BV algebra with $\Delta$ defined by (\ref{def-Laplace})
 and its definition requires the choice of a volume form on $M$.
The graded manifold  $T^*[-1]M$ is called a BV manifold. In general a BV manifolds is
 defined as a graded manifold ${\cal M}$ such that the space of functions $C^\infty ({\cal M})$ is equipped
  with the structure of a BV algebra.

 There also exists an alternative definition of BV algebra \cite{Getzler:1994yd}.
\begin{Def}\label{BValgebra-alternative}
A graded commutative algebra $V$ with an odd $\mathbb{R}$--linear map
 $$\Delta ~:~V \longrightarrow V~,$$
  which squares to zero $\Delta^2=0$ and satisfies
  \bea
  \Delta (vwz)&= & \Delta (vw) z + (-1)^{|v|} v \Delta (wz) + (-1)^{(|v|-1)|w|} w\Delta (vz) \nn \\
   &&- \Delta (v) wz -
   (-1)^{|v|} v \Delta(w) z - (-1)^{|v|+|w|} vw \Delta (z), \label{Delta-quadratic}
   \eea
  is called a BV algebra
 \end{Def}
  One can show that $\Delta$ with these properties gives rise to the bracket (\ref{def-bracket-Delta}) which satisfies
   all axioms of the definition \ref{Gerstenhaber-algebra}. The condition (\ref{Delta-quadratic}) is related to the fact that
    $\Delta$ should be a second order operator, square of the derivation in other words.
     Consider the functions $f(x)$, $g(x)$ and $h(x)$ of one variable and the second derivative
     $\frac{d^2}{dx^2}$ satisfies the following property
 $$ \frac{d^2 (fgh)}{dx^2}   +  \frac{d^2f}{dx^2}    gh + f\frac{d^2g}{dx^2}  h + fg \frac{d^2h}{dx^2}  =
  \frac{d^2(fg)}{dx^2} h + \frac{d^2( fh)}{dx^2} g + f \frac{d^2 (gh)}{dx^2} ~, $$
   which can be regarded as a definition of second derivative. Although one should keep in mind that any
    linear combination $\alpha \frac{d^2 }{d x^2} + \beta \frac{d}{dx}$ satisfies the above identity.
   Thus the property (\ref{Delta-quadratic})
    is just the graded generalization of the second order differential operator.
     In the example of $C^\infty (T^*[-1]M)$, the $\Delta$ as in (\ref{def-Laplace}) is indeed
     of second order.

     We collect some more details and curious observations on odd Fourier transform and some of its algebraic structures
      in  Appendices \ref{app-A} and \ref{app-B}.

 \subsection{Integration theory}\label{s-integration}

 So far we have discussed different algebraic aspects of graded manifolds $T[1]M$ and $T^*[-1]M$ which can be
  related by the odd Fourier transformation upon the choice of a volume form on $M$.  As we saw  $T^*[-1]M$
   is quite interesting algebraically since $C^\infty (T^*[-1]M)$ is equipped with the structure of a BV algebra.
    At the same time  $T[1]M$ has a very natural integration  theory which we will review below. Now our goal is
     to mix the algebraic aspects of $T^*[-1]M$ with the integration theory on $T[1]M$. We will do it again by means of
      the odd Fourier transform.

We start by reformulating the Stokes theorem in the language of the graded (super) manifolds. Before doing this let
 us review a few facts about standard submanifolds. A submanifold $C$ of $M$ can be described in algebraic language as
  follows. Consider the ideal ${\cal I}_C \subset C^\infty (M)$  of functions vanishing on $C$. The functions on submanifold
   $C$ can be described as quotient $C^\infty (C) = C^\infty (M) / {\cal I}_C$.  Locally we can choose coordinates $x^\mu$
    adapted to $C$ such that the submanifold $C$ is defined by the conditions $x^{p+1}=0, x^{p+2}=0~, ..., x^n=0$
     ($\dim C = p$ and $\dim M=n$) while the rest $x^1, x^2, .... , x^p$ may serve as coordinates for $C$. In this local description
      ${\cal I}_C$ is generated by  $x^{p+1}, x^{p+2}, ..., x^n$.  Indeed the submanifolds can be defined purely algebraically
        as ideals of  $C^\infty (M)$ with certain regularity condition which states that locally the ideals generated
 by $x^{p+1}, ... , x^n$. This construction has a straightforward generalization for the graded and
  super settings.  Let us illustrate this with a particular example which is relevant for our later discussion. $T[1]C$ is a graded submanifold of
   $T[1]M$ if $C$ is submanifold of $M$. In local coordinates $T[1]C$ is described by the conditions
   \bea
   x^{p+1}=0, ~~x^{p+2}=0~,~~ ...~, ~~x^n=0~,~~\theta^{p+1}=0~,~~ \theta^{p+2}=0~,~~ ...~, ~\theta^n=0~,
   \label{local-tangent}
   \eea
  thus $x^{p+1}, ... , x^n, \theta^{p+1}, ... ,\theta^n$ generate  the corresponding ideal ${\cal I}_{T[1]C}$.  The functions on the submanifold $C^\infty (T[1]C)$ are given by
    the quotient $C^\infty (T[1]M)/{\cal I}_{T[1]C}$. Moreover the above conditions define the morphism $i:T[1]C \rightarrow T[1]M$
     of the graded manifolds and thus we can talk about the pull back of functions from $T[1]M$ to $T[1]C$ as going to the quotient.
       Also we want to discuss another class of submanifolds, namely odd conormal bundle $N^*[-1]C$ as graded
      submanifold of $T^*[-1]M$. In local coordinate $N^*[-1]C$ is described by the conditions
   \bea
   x^{p+1}=0, ~~x^{p+2}=0~,~~ ...~, ~~x^n=0~,~~\psi_{1}=0~,~~ \psi_{2}=0~,~~ ...~, ~\psi_p=0~,\label{local-cotangent}
   \eea
 thus  $x^{p+1}, ... , x^n, \psi_1, ...  ,\psi_p$ generate the ideal ${\cal I}_{N^*[-1]C}$. Again the functions $C^\infty (N^*[-1]C)$ can be described as
     quotient $C^{\infty} (T^*[-1]M)/{\cal I}_{N^*[-1]C}$.  Moreover the above conditions define the morphism $j:N^*[-1]C \rightarrow T^*[-1]M$
     of the graded manifolds and thus we can talk about the pull back of functions from $T^*[-1]M$ to $N^*[-1]C$.

In previous subsections we have defined the odd Fourier transformation as map
$$ C^{\infty} (T[1]M)~~~~\overset{F}{\longrightarrow}~~~~C^{\infty}(T^*[-1]M)~,$$
 which does not map the graded commutative product on one side
  to the graded commutative product on the other side.  Using the odd Fourier transform we can relate
   the following integrals over different supermanifolds
  \bea
  \int\limits_{T[1]C} d^pxd^p\theta~ i^*\left ( f(x,\theta)\right )=(-1)^{(n-p)(n-p+1)/2}\int\limits_{N^*[-1]C} d^pxd^{n-p}\psi~\rho~ j^*\left (F[f](x,\psi)\right )~.
  \label{tangent-conormal-int}
  \eea
 Let us spend some time explaining this formula. On the left hand side we integrate the pull back of $f \in C^\infty (T[1]M)$ over
  $T[1]C$ with the canonical measure $d^pxd^p\theta$, where $d^p\theta = d\theta^p d\theta^{p-1} ... d\theta^1$.
   On the right hand side of (\ref{tangent-conormal-int}) we integrate the pull back of $F[f] \in C^\infty(T^*[-1]M)$ over
    $N^*[-1]C$.  The supermanifold  $N^*[-1]C$ has measure $d^px~ d^{n-p}\psi~ \rho$, where $d^{n-p}\psi = d\psi_{n} d\psi_{n-1}   ... d\psi_{p+1}$
     and we have to make sure that this measure is invariant under the change of coordinates which preserve $C$.
      Indeed this is easy to check. Let us take the adapted coordinates $x^\mu = (x^i, x^\alpha)$ such that $x^i$ ($i,j =1,2,..., p$) are the
       coordinates along $C$ and $x^\alpha$ ($\alpha, \beta, \gamma= p+1, ... , n$) are coordinates transverse to $C$. A generic
        change of coordinates has the form
         \bea
         \tilde{x}^i = \tilde{x}^i (x^j, x^\beta)~,~~~~~~~\tilde{x}^\alpha = \tilde{x}^\alpha (x^j, x^\beta)~,
        \eea
        if furthermore we want to consider the transformations preserving $C$ then the following conditions should be
         satisfied
         \bea
          \frac{\partial \tilde{x}^\alpha}{\partial x^i} (x^j, 0)=0~.
         \eea
     These conditions follow from the general transformation of differentials
         \bea
          d\tilde{x}^\alpha = \frac{\partial \tilde{x}^\alpha}{\partial x^i} (x^j, x^\gamma) d x^i +
           \frac{\partial \tilde{x}^\alpha}{\partial x^\beta} (x^j, x^\gamma)
           dx^\beta
         \eea
       and the Frobenius  theorem which states that $d\tilde{x}^\alpha$ should only go to $dx^\beta$ once restricted on $C$.
          In this case the adapted coordinates transform to adapted coordinates. On $N^*[-1]C$ we have the following transformations
         of  odd conormal coordinate $\psi_\alpha$
     \bea
      \tilde{\psi}_\alpha = \frac{\partial x^\beta}{\partial \tilde{x}^\alpha} (x^i, 0) \psi_\beta~.
     \eea
      Let us stress that $\psi_\alpha$ is a coordinate on $N^*[-1]C$ not a section, and the invariant
       object will be $\psi_\alpha dx^\alpha$.
    Under the above transformations restricted to $C$ we have the following property
    \bea
      d^p x~ d^{n-p} \psi ~\rho (x^i, 0) =   d^p \tilde{x} ~d^{n-p} \tilde{\psi} ~\tilde{\rho} (\tilde{x}^i, 0)~,
    \eea
     where, for the transformation of $\rho$ see the example \ref{integrationcotangent}.
       The formula  (\ref{tangent-conormal-int}) is very easy to prove in the local coordinates. The pull back of the functions
   on the left and right hand sides would correspond to imposing the conditions (\ref{local-tangent})
    and (\ref{local-cotangent}) respectively. The rest is just simple manipulations with the odd integrations
     and with the explicit form of the odd Fourier transform. Since all operations in (\ref{tangent-conormal-int}) are
      covariant, i.e. respects the appropriate gluing then the formula obviously is globally defined and is independent
       from the choice of the adapted coordinates.

 Let us recall two important corollaries of the  Stokes theorem for the differential forms. First corollary is that
  the integral  of exact form  over closed submanifold $C$ is zero and the second corollary is that
   the integral over closed form depends only on  homology class of $C$,
  \bea
   \int\limits_C d\omega =0~,~~~~~~ \int\limits_C \alpha = \int\limits_{\tilde{C}} \alpha~,~~d\alpha=0~,
  \eea
   where $\alpha$ and $\omega$ are differential forms, $C$ and $\tilde{C}$ are closed submanifolds which are
    in the same homology class.  These two statements can be easily rewritten in the graded language as
     follows
 \bea
  \int\limits_{T[1]C} d^pxd^p\theta~ D g=0~,\label{Stokes-th1}
 \eea
 \bea
  \int\limits_{T[1]C} d^pxd^p\theta~ f =   \int\limits_{T[1]\tilde{C}} d^pxd^p\theta~ f~,~~~~Df=0~,\label{Stokes-th2}
 \eea
  where we assume that we deal with the pull backs of $f,g \in C^\infty (T[1]M)$ to the submanifolds.

  Next we can combine the formula (\ref{tangent-conormal-int}) with the Stokes
    theorem (\ref{Stokes-th1}) and (\ref{Stokes-th2}). We will end up with the following properties to
     which we will refer as \emph{Ward-identities}
    \bea
    \int\limits_{N^*[-1]C} d^pxd^{n-p}\psi~\rho~ \Delta \tilde{g} =0~,\label{Ward1}
 \eea
 \bea
 \int\limits_{N^*[-1]C} d^pxd^{n-p}\psi~\rho~ \tilde{f} =
\int\limits_{N^*[-1]\tilde{C}} d^pxd^{n-p}\psi~\rho~ \tilde{f}~,~~~~~\Delta \tilde{f}=0~,\label{Ward2}
 \eea
  where $\tilde{f}, \tilde{g} \in C^\infty(T^*[-1]M)$ and the pull back of these function to $N^*[-1]C$ is assumed.
   One can think of these statements as a version of Stokes theorem for the cotangent bundle. This can be
   reformulated and  generalized further as  a general theory of integration over Lagrangian submanifold of
     odd symplectic supermanifold (graded manifold), for example see \cite{Schwarz:1992nx}.

  \subsection{Algebraic view on the integration}
\label{AVOTI}

Now we would like to combine the two facts about the graded cotangent bundle $T^*[-1]M$. From one side we have
 the BV-algebra structure on $C^\infty(T^*[-1]M)$, in particular we have the odd Lie bracket on the functions.
  From the other side we showed in the last subsection that  there exists an integration theory for $T^*[-1]M$ with an analog
   of the Stokes theorem.  Our goal is to combine the algebraic structure on $T^*[-1]M$ with the integration and argue that
    the integral can be understood as certain cocycle.

    Before discussing our main topic, let us remind the reader of some facts about the
     Chevalley-Eilenberg complex for the Lie
     algebras. Consider a Lie algebra $\FR{g}$ and define the space of $k$-chains $c_k$  as an element of $\wedge^k \FR{g}$.
 The space  $\wedge^k \FR{g}$ is  spanned by
     \bea
       c_k = T_1 \wedge T_2 \wedge ... \wedge T_k
     \eea
   and the boundary operator can be defined as follows
   \bea
    \partial (T_1 \wedge T_2 \wedge ... \wedge T_k) = \sum\limits_{1 \leq i < j \leq  k} (-1)^{i+j+1} [T_i, T_j] \wedge T_1 \wedge ... \wedge \hat{T}_i \wedge ... \wedge \hat{T}_j \wedge ...
     \wedge T_n~,
   \eea
    where $\hat{T}_i$ indicates the omission of the argument $T_i$ .
    Using the Jacobi identity one can easily prove that $\partial^2=0$. The dual object $k$-cochain $c^k$ is defined
     as multilinear map $c^k : \wedge^k \FR{g} \rightarrow \mathbb{R}$ such that coboundary operator $\delta$ is defined
      as follows
      \bea
        \delta c^k (T_1 \wedge T_2 \wedge ... \wedge T_k) = c^k \left ( \partial (T_1 \wedge T_2 \wedge ... \wedge T_k ) \right )
      \eea
       and $\delta^2=0$. This gives rise to the famous Chevalley-Eilenberg complex.  If $\delta c^k=0$ then we call $c^k$ a cocycle.
        If there exists $\tilde{c}^{k-1}$ such that $c^k = \delta \tilde{c}^{k-1}$ then we call $c^k$ a coboundary. The Lie algebra cohomology
         $H^k (\FR{g},\mathbb{R})$ consists the cocycles modulo coboundaries.
       In general we can also generalize it such
        that cochains take value in a $\FR{g}$-module. However this generalization is not relevant for the present discussion.

        Now let us consider the generalization of Chevalley-Eilenberg complex for the graded Lie algebras. Notice in the preceding paragraph we have defined the cochain as a mapping from $\wedge^k\FR{g}$ to numbers which is identified with   $\wedge^k\FR{g}^*$.
        However  $\wedge^k\FR{g}^*$ can also be thought of
          as $S (\FR{g}[1])$-the symmetric algebra of $\FR{g}[1]$ (see the definition \ref{def-symmalg-graded}).
           It is this formulation that allows for the most economical generalization to the graded case. Let
        a graded vector space $V$ equipped with Lie bracket $[~,~]$ of \emph{degree $0$}.  The cochains are
         defined as maps from graded symmetric algebra $S(V[1])$ to real numbers. More precisely, $k$-cochain is defined
          as multilinear map $c^k(v_1, v_2, ..., v_k)$ with the following symmetry properties
          \bea
           c^k(v_1,..., v_i, v_{i+1}, ..., v_k) = (-1)^{(|v_i|+1)(|v_{i+1}|+1)}   c^k(v_1,..., v_{i+1}, v_{i}, ..., v_k) ~.\label{sym_prop}
          \eea
       The coboundary operator $\delta$ is acting as follows
       \bea
        \delta c^k (v_1, ...v_{k+1})& = &\sum (-1)^{s_{ij}} c^k \big((-1)^{|v_i|}[v_i, v_j], v_1, ... , \hat{v}_i, ... , \hat{v}_j, ... , v_{k+1}\big)~,\nn\\
        s_{ij}&=&(|v_i|+1)(|v_1|+\cdots+|v_{i-1}|+i-1)\nn\\
        &&+(|v_j|+1)(|v_1|+\cdots+|v_{j-1}|+j-1)+(|v_i|+1)(|v_j|+1)~.\label{CE_diff}\eea
        The sign factor $s_{ij}$ is called the Kozul sign; it is incurred by moving $v_i,v_j$ to the very front. While the sign $(-1)^{|v_i|}[v_i, v_j]$ ensures that this quantity conforms to the symmetry property of (\ref{sym_prop}) when exchanging $v_i,v_j$. As before
         we use the same terminology, coboundaries and cocycles. The cohomology $H^k (V, \mathbb{R})$ is $k$-cocycles modulo
          $k$-coboundaries.

Now let us consider the case where the bracket is of \emph{degree 1}.  The corresponding cochains and coboundary operator
 can be defined using the parity reversion functor applied for the even Lie algebra.  Let us define $W=V[1]$ be graded vector
  space with Lie bracket of degree $1$.  Then  $k$-cochain is defined
          as multilinear map $c^k(w_1, w_2, ..., w_k)$ with the following symmetry properties
          \bea
           c^k(w_1,..., w_i, w_{i+1}, ..., w_k) = (-1)^{|w_i| |w_{i+1}|}   c^k(w_1,..., w_{i+1}, w_{i}, ..., w_k) ~.\label{1-bracket-sym_prop}
          \eea
       The coboundary operator $\delta$ is acting as follows
       \bea
        \delta c^k (w_1, ...w_{k+1})& = &\sum (-1)^{s_{ij}} c^k \big((-1)^{|w_i|-1}[w_i, w_j], w_1, ... , \hat{w}_i, ... , \hat{w}_j, ... , w_{k+1}\big),\nn\\
        s_{ij}&=& |w_i| (|w_1|+\cdots+|w_{i-1}|)\nn\\
        &&+|w_j|(|w_1|+\cdots+|w_{j-1}|)+ |w_i| |w_j|~.\label{1-bracket-CE_diff}\eea
     The formulas  (\ref{1-bracket-sym_prop})  and (\ref{1-bracket-CE_diff}) are obtained by the parity shift from the the case with the even
      bracket, i.e. from the formulas (\ref{sym_prop}) and (\ref{CE_diff}).  The cocycles, coboundaries and cohomology are defined
       as usual.

Now using the definition of cocycle for the odd Lie bracket let us state  the important theorem about the integration which is a simple
 consequence of Stokes theorem for the multivector fields (\ref{Ward1}) and (\ref{Ward2}).

\begin{Thm}\label{theorem-cocycle-BV}
   Consider a collection of functions  $f_1, f_2, ... , f_k \in C^\infty (T^*[-1]M)$ such that $\Delta f_i =0$ ($i=1,2,..., k$).
     Define the integral
 \bea
 c^k (f_1, f_2, ... , f_k; C) =  \int\limits_{N^*[-1]C} d^pxd^{n-p}\psi~\rho~ f_1 (x, \psi) ... f_k(x, \psi)~,\label{int-BV-cocycle}
 \eea
  where $C$ is closed submanifold of $M$. Then  $c^k (f_1, f_2, ... , f_k; C)$ is a cocycle with respect to the odd
  Lie algebra structure $( C^\infty (T^*[-1]M), \{~,~\})$
  $$ \delta c^k (f_1, f_2, ... , f_k; C)=0~.$$
   Moreover $c^k (f_1, f_2, ... , f_k; C)$ differs from $c^k (f_1, f_2, ... , f_k; \tilde{C})$ by a coboundary if $C$ is homologous to $\tilde{C}$, i.e.
   $$c^k (f_1, f_2, ... , f_k; C) - c^k (f_1, f_2, ... , f_k; \tilde{C}) = \delta \tilde{c}^{k-1}~,$$
    where $\tilde{c}^{k-1}$ is some $(k-1)$-cochain.
\end{Thm}

This theorem is based on the observation by A.~Schwarz in \cite{Schwarz:1999vn}.  Let us now present the proof of this theorem.
$C^\infty (T^*[-1]M)$ is a graded vector space with odd Lie bracket $\{~,~\}$ defined in (\ref{def-PBoncot}), and the functions with
 $\Delta f =0$ correspond to a Lie subalgebra of $C^\infty (T^*[-1]M)$. The integral (\ref{int-BV-cocycle}) defines a $k$-cochain
  for odd Lie algebra with the correct symmetry properties
\bea
 c^k (f_1,  ... , f_i, f_{i+1}, ...  , f_k; C) = (-1)^{|f_i| |f_{i+1}|} c^k (f_1,  ... , f_{i+1}, f_{i}, ...  , f_k; C) ~,\nn
 \eea
  which follows from the graded commutativity of $C^\infty (T^*[-1]M)$.
 Then the property (\ref{Stokes-th1}) implies the following
 \bea
 0=\int\limits_{N^*[-1]C} d^pxd^{n-p}\psi~\rho~ \Delta(f_1 (x, \psi) ... f_k(x, \psi))~.\label{int-BV-cocycle-extra}
 \eea
Using the property of $\Delta$ given in (\ref{Delta-defPB}) many times, we obtain the following formula
\bea\label{lale03030}
 \Delta (f_1 f_2 ...  f_{k})&=&\sum_{i<j}\textrm(-1)^{s_{ij}}(-1)^{|f_i|}~  \{f_i,f_j\}~ f_1~ ...
~\widehat{f_i} ~ ... ~ \widehat{f_j} ~ ... ~ f_k ~,\nn\\
s_{ij}&=&(-1)^{(|f_1|+\cdots+|f_{i-1}|)|f_i|+(|f_1|+\cdots
+|f_{j-1}|)|f_j|-|f_i||f_j|}~,\label{Delta-to-sum}\eea
 where we have used $\Delta f_i =0$.
  Combining (\ref{int-BV-cocycle-extra}) and (\ref{Delta-to-sum}) we obtain that $c^k$ defined in (\ref{int-BV-cocycle}) is a cocycle, i.e.
  \bea
  \delta c^{k}(f_1, ... , f_{k+1}; C)=  - \int\limits_{N^*[-1]C} d^pxd^{n-p}\psi~\rho~ \Delta(f_1 (x, \psi) ... f_k(x, \psi))=0~,
  \eea
   where we use the definition for coboundary operator given in (\ref{1-bracket-CE_diff}).

 Next we have to show that the cocycle (\ref{int-BV-cocycle}) changes by a coboundary when we deform $C$ continuously.
  We start by looking at the infinitesimal change of $C$.  Recall that the indices $i,j$ are along $C$, while $\alpha,\beta$ are transverse to $C$, and $N^*[-1]C$ is locally given by $x^{\alpha}=0$ and $\psi_i=0$.
 We parameterize the infinitesimal deformation by
\bea
\delta_C x^{\alpha}=\epsilon^{\alpha}(x^i)~,~~~~\delta_C \psi_i=-\partial_i\epsilon^{\alpha}(x^i)\psi_{\alpha}~,\nn
\eea
 where $\epsilon$'s parametrize the deformation.
Thus a function $f \in C^\infty (T^*[-1]M)$ it changes as follows
\bea
 \delta_C f(x,\psi)\big|_{N^*[-1]C}=\epsilon^{\alpha}\partial_{\alpha}f-\partial_i\epsilon^{\alpha}(x^i)\psi_{\alpha}\partial_{\psi_{j}}f\Big|_{N^*[-1]C}
=-\{\epsilon^{\alpha}(x^i)\psi_{\alpha},f\}\big|_{N^*[-1]C}~.\nn\eea
Using $\Delta$ the bracket can be rewritten as
\bea
\delta_C f(x,\psi)\big|_{N^*[-1]C}=\Delta(\epsilon^{\alpha}(x^i)\psi_{\alpha}f)+\epsilon^{\alpha}(x^i)\psi_{\alpha}\Delta(f)\big|_{N^*[-1]C}~.\nn\eea
The first term vanishes under the integral. Specializing to infinitesimal deformation of $f_1\cdots f_k$,
we have
\bea &&\delta_Cc^k (f_1,..., f_k; C)=\delta\tilde c^{k-1} (f_1, ... , f_k)~,\nn\\
\textrm{where}:&&\tilde c^{k-1} (f_1,\cdots f_{k-1}; C)= - \int\limits_{N^*[-1]C}d^pxd^{n-p}\psi~\rho~
 \epsilon^{\alpha}(x^i)\psi_{\alpha}~f_1\cdots f_{k-1}~.\nn\eea
This shows that under an infinitesimal change of $C$, $c^k$ changes by a coboundary. For finite deformations of $C$, we can parameterize the deformation as a one-parameter family $C(t)$. Thus we have the identity
\bea \frac{d}{dt}c^k (f_1, ... , f_k; C(t))=\delta \tilde c^{k-1}(f_1, ... , f_k; C(t)) \nn\eea
for every $t$, integrating both sides we finally arrive at the the formula for the finite change of $C$
\bea c^k (f_1, ... , f_k; C(1))-c^k (f_1, ... , f_k; C(0)) =\delta\int\limits_0^1dt~  \tilde c^{k-1}_{C(t)}~.\nn\eea
This concludes the proof of Theorem \ref{theorem-cocycle}. $\blacksquare$

Now let us perform the integral (\ref{int-BV-cocycle}) explicitly.  Assume that the functions $f_i$ are of fixed degree
and we will use the same notation for the corresponding multivector  $f_i \in \Gamma (\wedge^\bullet TM)$. After
 pulling back the $f$'s and performing the odd integration the expression (\ref{int-BV-cocycle}) becomes
  \bea
   c^k (f_1, .... , f_k; C) = \int\limits_C i_{f_1} i_{f_2} ... i_{f_k} {\rm vol}~,\label{vector-diff-form}
  \eea
   where $i_f$ stands for the contraction of multivector with a differential form. Here all vector fields are assumed to
    be divergenceless.
   If $n -p \neq |f_1|+ ... + |f_k|$ then this
    integral is identically zero, otherwise it may be non-zero and it gives rise to the cocycle on
    $\Gamma (\wedge^\bullet TM)$ equipped with the Schouten bracket.

  It is easy to construct the examples of cocycles when we restrict our attention to the vector fields only.

 \begin{Exa}\label{vector-fields-cocycle}
  Let us illustrate the formula (\ref{vector-diff-form}) with a concrete set of examples. Let us pick up the volume form
   $\rm{vol}$ on $M$ and the collection of vector fields $v_i \in \Gamma (TM)$ which preserve this volume form, i.e.
    ${\cal L}_{v_i} \rm{vol}=0$, where ${\cal L}_{v_i}$ is a Lie derivative with respect to ${\mathbf v}_i$.
     Such vector fields form
     the closed Lie algebra and they are automatically divergenceless. If $n-p=k$ then the integral
 $$  c^k (v_1, .... , v_k; C) = \int\limits_C i_{v_1} i_{v_2} ... i_{v_k} {\rm vol}~,\label{vector-diff-form-example}$$
  gives rise to the cocycle of the Lie algebra of the vector field preserving a given volume form.
   Such situation is realized in many examples. For instance, consider Lie group ${\mathbf G}$ with the left invariant
    vector fields and left invariant volume form, Hamiltonian vector fields on Poisson manifold with unimodular Poisson
     structure, the Hamiltonian vector fields on symplectic manifold etc.
    However one should keep in mind that the corresponding cocycles are relatively trivial in some sense.
 \end{Exa}

 The main moral of Theorem \ref{theorem-cocycle-BV} is that the BV integral gives rise to cocycle with specific dependence from $C$
  and indeed it can be taken as a defining property of those integrals. In infinite dimensional setting when the integral is not defined
   at all then the statement of this Theorem can be regarded as definition. We will briefly comment on it in section \ref{QFT-summary}.

\section{Perturbation theory}
\label{perturbative}

The goal of this section is to give an introduction to the perturbative expansion of finite dimensional integrals
 and discuss the element of graph theory which are relevant for further discussion.
  We will also state the Kontsevich theorem.
  This section should be
  regarded as a technical preparation for the next section.

Gaussian integrals, though quite elementary, are really at the very core of perturbation theory. We quickly go over the Gaussian integrals on
 $\BB{R}^n$, focusing on how to organize the calculation in terms of Feynman diagrams. We will be quite brief in our discussion of
  the perturbative expansion. For mathematically minded reader we recommend  two nice short introductions
   to perturbation theory, \cite{Polyak:2004nn} and \cite{Sawon:2005qt}.

\subsection{Integrals in $\mathbb{R}^n$-Gaussian Integrals and Feynman Diagrams}\label{IiRnGIaFD}

 Let us discuss how to calculate the specific integrals on $\mathbb{R}^n$. We are interested in the combinatorial and algebraic
  way of the calculation of the integrals. We recall that the standard one-dimensional Gaussian integral is given by
   the following formula
 \bea
    \int\limits_{-\infty}^{\infty} dx~ e^{-\frac{1}{2} \alpha x^2} = \sqrt{ \frac{2\pi}{\alpha}}~.
 \eea
 The corresponding generalization for $\mathbb{R}^n$ is given by
 \bea
 \int\limits_{\mathbb{R}^n} d^nx~e^{-\frac{\alpha}{2} Q_{\mu\nu} x^\mu x^\nu} = \left ( \frac{2\pi}{\alpha}\right )^{n/2} \frac{1}{\sqrt{\det Q}} \equiv Z[0]~,\label{R-n-gauss-result}
 \eea
  where $Q_{\mu\nu}$ is a symmetric and positive matrix.  Next we introduce the generating function of $n$-variables $J_1, J_2, ...., J_{n}$
   as follows
 \bea
  Z[J]  = \int\limits_{\mathbb{R}^n} d^nx~e^{-\frac{\alpha}{2} Q_{\mu\nu} x^\mu x^\nu + J_\mu x^\mu} = Z[0] ~e^{\frac{1}{2\alpha} Q^{\mu\nu} J_\mu J_\nu}~,
 \eea
  where $Q^{\mu\nu} Q_{\nu\lambda} = \delta^\mu_\lambda$.  Let us introduce the integral
 \bea
 \langle x^\sigma x^\lambda \rangle = \frac{1}{Z[0]} \int\limits_{\mathbb{R}^n} d^nx~x^\sigma x^\lambda~ e^{-\frac{\alpha}{2} Q_{\mu\nu} x^\mu x^\nu} = \frac{1}{Z[0]} \frac{\partial^2}{\partial J_\sigma \partial J_\lambda} Z[J] |_{J=0} =\frac{1}{\alpha} Q^{\sigma\lambda}~,
 \eea
  which is straightforward to calculate. Next we would like to discuss the following integrals
 \bea
 \langle x^{\mu_1} x^{\mu_2} ... ~x^{\mu_{2n}}  \rangle = \frac{1}{Z[0]} \int\limits_{\mathbb{R}^n} d^nx~x^{\mu_1} x^{\mu_2} ... ~ x^{\mu_{2n}}~ e^{-\frac{\alpha}{2} Q_{\mu\nu} x^\mu x^\nu}~,\label{many-x-int}
 \eea
  first if we have an odd number of $x$'s under the integral then it is identically zero due to symmetry properties.  The integral
   (\ref{many-x-int}) can be calculated  by taking the derivatives
   \bea
    \langle x^{\mu_1} x^{\mu_2} ... ~x^{\mu_{2n}}  \rangle = \frac{\partial^{2n}}{\partial J_{\mu_1} \partial J_{\mu_2} ... \partial J_{\mu_{2n}}}
      e^{\frac{1}{2\alpha} Q^{\mu\nu} J_\mu J_\nu} |_{J=0}~.
    \eea
  Performing the derivatives explicitly we arrive at the statement which is known in physics literature as the Wick theorem,
  \bea
       \langle x^{\mu_1} x^{\mu_2} ... ~x^{\mu_{2n}}  \rangle = \frac{1}{2^n n! \alpha^n} \sum\limits_{P} ~Q^{\mu_{P_1}\mu_{P_2}}~
       ...~ ~Q^{\mu_{P_{2n-1}} \mu_{P_{2n}}}~,\label{Wick-theorem-expl}
  \eea
 where we sum over all permutations $P$ of the indices $\mu_1, \mu_2, ... , \mu_{2n}$.
   In general we are interested in the following integrals
 \bea
  \langle V_{n_1} (x) .... V_{n_k} (x) \rangle =   \frac{1}{Z[0]} \int\limits_{\mathbb{R}^n} d^nx~ V_{n_1}(x) ... ~ V_{n_k}(x)~ e^{-\frac{\alpha}{2} Q_{\mu\nu} x^\mu x^\nu}~,\label{general-int-calc}
 \eea
 where $V$'s are monomials in $x$ of the fixed degree
 \bea
  V_{n} (x)  = \frac{1}{n!} V_{\mu_1 ... \mu_n} x^{\mu_1} ... x^{\mu_n}
 \eea
  and let us assume for the moment that they are all different. Using the Wick theorem (\ref{Wick-theorem-expl}) we can calculate
  \bea
    \langle V_{n_1} (x) .... V_{n_k} (x) \rangle = \frac{1}{n_1! ... n_k!} V_{\mu_1 ...\mu_{n_1}} ... ~V_{\nu_1 ... \nu_{n_k}} \langle x^{\mu_1} ~...~~ x^{\nu_{n_k}} \rangle~,
  \eea
   where
   we have to contract $V$'s with $Q$'s in all possible ways. The ways to contract $V$'s with $Q$ can be depicted using the graph
     $\Gamma$  (Feynman diagram), where $V_{\mu_1...\mu_n}$ is $n$-valent vertex and $Q^{\mu\nu}$ corresponds to the edges.
      Following physics terminology we will call $V$ a vertex and $Q^{\mu\nu}$ a propagator.
 Thus  in perturbation theory, the integrals (\ref{general-int-calc}) are effectively organized as Feynman diagrams, and physicists have developed effective mnemonic rules to keep track of the combinatorics. Amongst these rules the symmetry factors $|{\rm Aut}~ \Gamma|$ of a Feynman diagram $\Gamma$ is the most important. Due to the close relationship between Feynman diagrams and the Kontsevich theorem about the graph complex, we spend some time to go through the examples of the Feynman diagrams. The general formula for the integral
   (\ref{general-int-calc})
  looks as follows
   \bea
       \langle V_{n_1} (x) .... V_{n_k} (x) \rangle = \frac{1}{\alpha^{(n_1+... +n_k)/2}}\sum\limits_{\Gamma} \frac{1}{|{\rm Aut}~ \Gamma|} W(\Gamma)~, \label{pert-gener-form}
\eea
 where $|{\rm Aut}~\Gamma|$ is the symmetry factor for a diagram $\Gamma$ and $W(\Gamma)$ is the contraction of $V$'s with $Q$'s
  according to the $\Gamma$.

It turns out the weight of a diagram is determined by the inverse of its symmetry factor $|{\rm Aut~\Gamma}|$.
 The symmetry factor consists of three parts $|{\rm Aut~\Gamma}| = (\#P) (\# V) (\# L)$,
\begin{itemize}
\item $\#P$ the symmetry factor of edges, if there are $p$ edges running between a pair of vertices we include a factor of $p!$.
\item $\#V$ the symmetry factor of vertices, defined as the cardinality of the subgroup of $s_k$ that preserves the graph (disregarding the orientation of edges)
\item $\#L$ in case loops (edges starting and ending on the same vertex) are allowed, a factor of 2 for each such loop
\end{itemize}

Now we present the examples which clarify and explain the formula (\ref{pert-gener-form})

 \begin{Exa}\label{example-one-dim-pert}
  The above logic is equally applicable for one-dimensional integrals.
 The following integral is  very easy to calculate
 \bea
   \int\limits_{-\infty}^{\infty} dx~ x^{2p}e^{-\frac{1}{2} \alpha x^2} = \sqrt{ \frac{2\pi}{\alpha}}\frac{(2p-1)!!}{\alpha^p} = Z[0] \frac{(2p-1)!!}{\alpha^p}~.\label{direct_calc}
 \eea
Thus we have the explicit  integrals
\bea
&&\frac{1}{2} \frac{1}{(3!)^2} \langle   x^3 x^3 \rangle = \frac{1}{Z[0]}
 \int_{\BB{R}}dx~ \frac1{2!}(\frac1{3!}x^3)(\frac1{3!}x^3) e^{-\frac{1}{2} \alpha x^2}= \frac{5!!}{2!3!3!\alpha^3}~,\label{37uus333}\\
&&\frac{1}{3! 5!} \langle x^3 x^5\rangle = \frac{1}{Z[0]}
 \int_{\BB{R}}dx~ (\frac1{3!}x^3)(\frac1{5!}x^5) e^{-\frac{1}{2} \alpha x^2}=  \frac{7!!}{5!3!\alpha^4}~.\label{2j2j333} \eea
In the integral (\ref{37uus333}) we have included a conventional factor of $\frac{1}{2}$ because we have two identical vertices.
 Now let us recalculate these integrals using the Wick theorem (\ref{Wick-theorem-expl}) and related combinatorics.
 Taking $x^3$ $x^5$ as 3 and 5-point vertices respectively, the Gaussian integral is, according to Wick's contraction rule, summing over all possible ways of connecting all legs of the two vertices together.

The first line in Figure \ref{Gauss_feyn_fig} representes the Feynman diagrams to connect $x^3$ to $x^3$
\bea
\frac{1}{2} \frac{1}{(3!)^2} \langle   x^3 x^3 \rangle = \frac{1}{2} \frac{1}{(3!)^2\alpha^3} (6 + 9) ~,
\eea
 where $6$ is the number of ways of contractions according to the first diagram and $9$ is the number of ways of contraction according to the second
  diagram. Altogether there is 15 ways of contracting $x^3$ with $x^3$, and the number $15$ may also be calculated as
   follows
   $$\frac{1}{3!} \left (\begin{array}{c} 6\\ 2~2~2 \end{array}\right)=15~.$$
   Furthermore we have
   \bea
   \frac{1}{2} \frac{1}{(3!)^2\alpha^3} (6 + 9) =  \frac{1}{\alpha^{3}}\left (\frac{1}{3!2}+\frac18\right )=\frac{5}{24\alpha^3}~,
   \eea
 which is exactly the same as the integral (\ref{37uus333}). Indeed  $3!2$ and $8$ are the symmetry factors  for the diagrams in the first line
 of Figure \ref{Gauss_feyn_fig}.  Thus we are in agreement with the general prescription given by the formula (\ref{pert-gener-form}). \\

\begin{figure}[h]
\begin{center}
\includegraphics[width=2.2in]{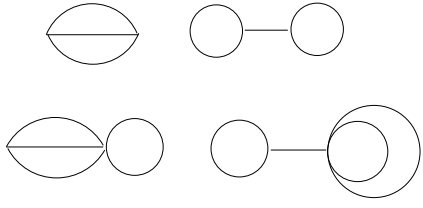}
\caption{Symmetry factor for the graphs are $12,8,12,16$ respectively}\label{Gauss_feyn_fig}
\end{center}
\end{figure}

Now let us calculate the integral with $x^3$ and $x^5$. Altogether there are 105 ways of contracting $x^3$ and $x^5$
  $$\frac{1}{4!} \left (\begin{array}{c} 8\\ 2~2~2~2 \end{array}\right)=105~.$$
Thus we have
\bea
\frac{1}{3! 5!} \langle x^3 x^5 \rangle = \frac{1}{3! 5! \alpha^4} (60 + 45)~,
\eea
 where $60$ corresponds to the first diagram in the second line and $45$ to the second diagram in the second line of Figure
   \ref{Gauss_feyn_fig}. Next we get
\bea \frac{1}{3! 5! \alpha^4} (60 + 45)= \frac{1}{\alpha^{4}}\left (\frac{1}{3!2}+\frac1{16}\right )=\frac{7}{48\alpha^4}~,\nn\eea
 which is exactly the same as the integral (\ref{2j2j333}). Indeed $3!2$ and $16$ are  symmetry factors for the diagrams in
  the second line  of  Figure \ref{Gauss_feyn_fig}. Again we have agreement with the prescription of  (\ref{pert-gener-form}).
\end{Exa}

\begin{Exa}\label{example-n-dim-pert}
  Next consider the $\mathbb{R}^n$-analogs of the previous example. Let us first calculate the following integral
    \bea
  \frac{1}{2} \langle V_3(x) V_3(x) \rangle \equiv \frac{1}{Z[0]} \int\limits_{\mathbb{R}^n} d^nx~  \frac{1}{2}\frac{1}{3!} V_{\mu_1 \mu_2 \mu_3} x^{\mu_1} x^{\mu_2} x^{\mu_3}
  \frac{1}{3!} V_{\mu_4 \mu_5 \mu_6} x^{\mu_4} x^{\mu_5} x^{\mu_6}~ e^{-\frac{\alpha}{2} Q_{\mu\nu} x^\mu x^\nu}~.
 \eea
  Using the Wick theorem and the counting from the previous example we get
 \bea
 \frac{1}{2}\langle V_3 (x) V_3(x) \rangle =  \frac{1}{2 (3!)^2}  \frac{1}{\alpha^3}\left ( 6 W (\Gamma_1) + 9 W (\Gamma_2) \right )~,
 \eea
 where now we have non-trivial factors
 \bea
  W(\Gamma_1) = V_{\mu_1\mu_2 \mu_3} V_{\mu_4 \mu_5 \mu_6} Q^{\mu_1 \mu_4} Q^{\mu_2 \mu_5} Q^{\mu_3 \mu_6}~,\nn
 \eea
 and
 \bea
  W(\Gamma_2) = V_{\mu_1\mu_2 \mu_3} V_{\mu_4 \mu_5 \mu_6} Q^{\mu_1 \mu_2} Q^{\mu_4 \mu_5} Q^{\mu_3 \mu_6}~.\nn
 \eea
  Moreover we can rewrite it as
  \bea
 \frac{1}{2}\langle V_3 (x) V_3(x) \rangle =
  \frac{1}{\alpha^3} \left (\frac{1}{12} W(\Gamma_1) + \frac{1}{8} W(\Gamma_2) \right )~,
 \eea
 where $12$ and $8$ are the symmetry factors for the diagrams in first line of Figure \ref{Gauss_feyn_fig}.
  Analogously we can calculate
 \bea
  \langle V_3(x) V_5(x) \rangle \equiv \frac{1}{Z[0]}\int\limits_{\mathbb{R}^n} d^nx~  \frac{1}{3!} V_{\mu_1 \mu_2 \mu_3} x^{\mu_1} x^{\mu_2} x^{\mu_3}
  \frac{1}{5!} V_{\mu_4 \mu_5 \mu_6 \mu_7 \mu_8} x^{\mu_4} x^{\mu_5} x^{\mu_6} x^{\mu_7} x^{\mu_8}~ e^{-\frac{\alpha}{2} Q_{\mu\nu} x^\mu x^\nu}~,
 \eea
  which after applying the Wick theorem gives us
 \bea
   \langle V_3(x) V_5(x) \rangle =  \frac{1}{3! 5!} \frac{1}{\alpha^4} \left ( 60 W(\Gamma_3) +  45 W(\Gamma_4) \right)~,
 \eea
where
\bea
W(\Gamma_3) =  V_{\mu_1\mu_2 \mu_3} V_{\mu_4 \mu_5 \mu_6 \mu_7 \mu_8} Q^{\mu_1 \mu_4} Q^{\mu_2 \mu_5} Q^{\mu_3 \mu_6}
 Q^{\mu_7 \mu_8}\nn
\eea
  and
  \bea
  W(\Gamma_4) = V_{\mu_1\mu_2 \mu_3} V_{\mu_4 \mu_5 \mu_6\mu_7 \mu_8} Q^{\mu_1 \mu_2} Q^{\mu_4 \mu_5} Q^{\mu_6 \mu_7}
   Q^{\mu_3 \mu_8}~.\nn
 \eea
This can be rewritten as
   \bea
   \langle V_3(x) V_5(x) \rangle =   \frac{1}{\alpha^4} \left ( \frac{1}{12} W(\Gamma_3) +  \frac{1}{16} W(\Gamma_4) \right)~,
 \eea
  where $12$ and $16$ are the symmetry factors for the diagrams in the second line of  Figure \ref{Gauss_feyn_fig}.
\end{Exa}

In these two examples we illustrated the prescription given by the formula (\ref{pert-gener-form}). If there are $p$ identical monomials
 $V$  then we have to put the additional $1/p!$ on the left hand side for this formula to work.

So far we have assumed that the matrix $Q_{\mu\nu}$ is positive definite and that all integrals, considered so far, do converge.
However we can treat the integrals formally and drop the condition of matrix  $Q_{\mu\nu}$ being positive definite.
Thus all our manipulations with the integrals  can be merely taken as a book keeping device for Wick contractions.
 In all following discussion we treat the integrals formally and will not ask if the integral converges at all.
  Of course, we have to keep in mind that many formal integrals can be understood less formally through an analytic continuation
   as convergent integrals. This can be an important point if we try to go beyond the perturbative expansions.

\subsection{Integrals in $\bigoplus\limits_{i=1}^N \mathbb{R}^{2n}$}\label{IiRN}

As a preparation for section \ref{AUBVToaL} we give some details  of the perturbative expansion of
the generalization of integrals from the previous subsection.
Let $\Omega_{\mu\nu}$ be the constant symplectic form on $\BB{R}^{2n}$ and let $t_{ij}= - t_{ji}$ $i,j=1\cdots N$ be
 an antisymmetric non-degenerate matrix and $t^{ij}$ its inverse. Let
 us define the formal Gaussian integral  over $N$ copies of $\BB{R}^{2n}$,  where each copy is labeled by a subscript.
  The integral is defined as follows
\bea
Z[0] =  \int\limits_{\bigoplus\limits_{i=1}^N \mathbb{R}^{2n}} d^{2n} x_1 ... d^{2n}x_N ~e^{-\frac{1}{2}\sum\limits_{i,j}t^{ij} \Omega_{\mu\nu} x_i^\mu x_j^\nu} = (2\pi)^{nN}(\det t)^{-n}(\det \Omega)^{-N/2}~,
\eea
 where we used (\ref{R-n-gauss-result}) which is now understood as formal expression.
From now on, we use the summation convention that any repeated Greek indices are summed while repeated Latin indices are not summed unless explicitly indicated. In general we are interested in the following integrals
\bea
 \int\limits_{\bigoplus\limits_{i=1}^N \mathbb{R}^{2n}} d^{2n} x_1 ... d^{2n}x_N ~f_1 (x_1) f_2(x_2) ... f_N (x_N)
  ~e^{-\frac{1}{2}\sum\limits_{i,j}t^{ij} \Omega_{\mu\nu} x_i^\mu x_j^\nu}~,\label{Gauss_feyn_high_d1}\eea
  where $f_i (x_i)$ are polynomials which depend on the coordinate of a copy of $\BB{R}^{2n}$.  Obviously it is enough to
   consider only monomials and we assume the following normalization
   \bea
   f(x) = \frac{1}{n!} ~f_{\mu_1 \mu_2 ... \mu_n}~ x^{\mu_1} x^{\mu_2} ... x^{\mu_n}~.\label{normalization-pol}
   \eea
   For physics minded readers we can comment that the above integral can be thought of as
 discrete field theory whose source manifold is $N$ points labeled by $i=1\cdots N$ and target is $\BB{R}^{2n}$, and whose 'fields' are $x_i$.

One can work out these integrals by usual methods,
\bea
&& \langle f_1(x_1) f_2(x_2) ... f_N(x_N) \rangle  \nn \\
 &\equiv  &\frac{1}{Z[0]} \int\limits d^{2n} x_1 ... d^{2n}x_N ~f_1 (x_1) f_2(x_2) ... f_N (x_N)
  ~e^{-\frac{1}{2}\sum\limits_{i,j}t^{ij} \Omega_{\mu\nu} x_i^\mu x_j^\nu}\nn\\
  &=&
\frac{1}{Z[0]} \int d^{2n} x_1 ... d^{2n}x_N ~f_1 (\partial_{J^1}) f_2(\partial_{J^2}) ... f_N (\partial_{J^N})
  ~e^{-\frac{1}{2}\sum\limits_{i,j}t^{ij} \Omega_{\mu\nu} x_i^\mu x_j^\nu+\sum\limits_i J^i_{\mu}x^{\mu}_i}\Big|_{J=0}\nn\\
&=&
\frac{1}{Z[0]}~f_1 (\partial_{J_1}) f_2(\partial_{J_2}) ... f_N (\partial_{J_N})
  ~e^{\sum\limits_{i,j}\frac12t_{ij}(\Omega^{-1})^{\mu\nu}J^i_{\mu}J^j_{\nu}}\Big|_{J=0}~, \label{Gauss_feyn_high_d}\eea
where  we apply the Wick theorem. Thus we see clearly that the answer can be represented as Feynman diagrams (graphs), for which the Taylor coefficient of $f_i$ serves as vertices and $t\Omega^{-1}$ serves as propagators.
 In this expansion the edges starting and ending at the same vertex are forbidden due to symmetry properties.
This is pretty straightforward to work out, but the reader may use formula (\ref{Gauss_feyn_high_d}) to understand where does the symmetry factor $\#P$ come from. There is no $\#V$ factor for this case since all vertices are distinct. In next section we will consider this perturbation
 theory further in the context of BV formalism and Kontsevich theorem.

 \subsection{Bits of graph theory}
 \label{BOGT}

 In previous subsections we encountered the graphs which depict the rules of contracting the indices of vertexes with propagators.
  In physics such graphs are called Feynman diagrams. In this subsection we would like to state very briefly some relevant
   notions of the graph theory.

By a graph we understand  a finite 1-dimensional CW complex.\footnote{We recall that the CW complex is a topological space whose building blocks are cells (space homeomorphic to $\BB{R}^{n}$ for some $n$), with the stipulation that each point in the CW complex is in the interior of one unique cell and the boundary of a cell is the union of cells of lower dimension. For example, the minimal cell decomposition of a sphere consists of one 0-cell, the north pole, and one 2-cell, containing the rest of the sphere.}
In simpler words, graph is collection of points (vertices) and lines (edges) with lines connecting the points. The reader may see the examples of
 the graphs on Figure \ref{Gauss_feyn_fig}.  We  consider only closed graphs, i.e. without external legs.  The graphs depicted on
  Figure \ref{Gauss_feyn_fig} are so called \emph{free graphs}. If we number the vertices by $1,2,...$ then we will call such graph
  \emph{labelled graph}.   For every labelled graph we can construct the adjacency matrix $\pi$ where matrix element   $\pi_{ij}$
    is equal to the number of edges between vertices $i$ and $j$.  The adjacency matrix $\pi$ can be used to describe the symmetries
     of the graph. If we relabel the vertices for a given graph then the adjacent matrix transforms by the similarity transformations
      $\tilde{\pi} = {\cal P} \pi {\cal P}^t$, where ${\cal P}$ is permutation matrix defined such ${\cal P}_{ij}=1$ if $i$ goes to $j$ under
       the permutation and otherwise zero.
       The permutation ${\cal P}$ is called a symmetry if the following satisfied $\pi = {\cal P} \pi {\cal P}^t$ (or in other words
       ${\cal P}$ and $\pi$ commute).
       Such ${\cal P}$'s give rise to the symmetry group of the graph
         and  the order of this group gives us $\#V$ which we defined previously as the symmetry factor of vertices.

An important structure of  graph is its
\emph{orientation}. The orientation is given by
\begin{itemize}
\item ordering of all the
vertices,
\item orienting of all the edges.
\end{itemize}
If one graph with $n$ vertices of a given orientation can be turned into another after a permutation $\sigma\in s_n$ of vertices ($s_n$ is the symmetry group of $n$-elements) and flipping of $k$ edges, then we say they are equal   orientation
 if $(-1)^k\textrm{sgn}(\sigma)=1$ and otherwise opposite orientation.
   We can talk about the \emph{oriented labelled graph}, which is the labelled graph with fixed order of vertices and oriented edges.
 We can introduce the \emph{equivalence classes of graphs} under the following relation
 $$(\Gamma, {\rm orientation})~~\sim~~(-\Gamma, {\rm opposite~ orientation})~.$$
  Thus every free graph comes with two orientations and we can now multiply the graphs (equivalence classes of graphs) by numbers
   and sum them formally.
 Thus the underlying module of the graph complex is generated by all the \emph{equivalence class of graphs}, and the grading of the complex is of course the number of vertices in a graph. Quite often to represent the equivalence class of graph we will use the concrete oriented
  labelled graph.

We remark here that there
are other orientation schemes that are equivalent to the current
one. For example, the orienting of all the edges can be replaced
by the ordering of all the legs from all vertices. Another
convenient scheme is to order the incident legs for each vertex
and order all the even valent vertices. We refer the reader to
section 2.3.1 of  the work \cite{ConantVogtmann} for the full discussion of  orientation.

The graph complex comes with a differential, which acts on the graph by shrinking  one edge and combining the two vertices connected by the edge. The main subtlety is to define the sign of such an operation wisely so as to make the differential nilpotent.
\begin{figure}[h]
\begin{center}
\psfrag{i}{\scriptsize{$i$}}\psfrag{j}{\scriptsize{$j$}}
\includegraphics[width=2in]{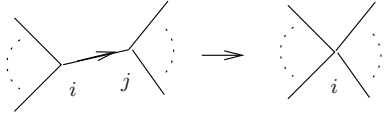}
\caption{The sign factor associated with this operator is $(-1)^j$}\label{graph_diff_fig}
\end{center}
\end{figure}
More precisely, when combining the vertices $i,j$ ($i<j$) with an edge from $i$ to $j$, we form a new vertex named $i$ inheriting all edges landing on $i,j$ (except the shrunk one of course), all the vertices labeled after $j$ move up one notch, and the sign factor associated to this procedure is $(-1)^j$. Note that in \cite{2010arXiv1006.1240Q}, the vertices are labeled starting from 0 instead of 1, so the sign factor there is $(-1)^{j+1}$. In case the edge runs from $j$ to $i$ one gets an extra $-$ sign.

\begin{figure}[h]
\begin{center}
\psfrag{del}{\scriptsize{$\partial_{Gph}$}}\psfrag{=6}{\scriptsize{$=6\times$}}\psfrag{=2}{\scriptsize{$=2\times$}}
\psfrag{G1}{\scriptsize{$\Gamma_1$}}\psfrag{G2}{\scriptsize{$\Gamma_2$}}\psfrag{G8}{\scriptsize{$\Gamma_8$}}
\includegraphics[width=2.2in]{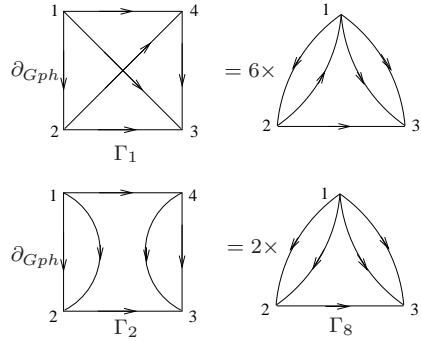}
\caption{Graph differential}\label{ex_box_diff}
\end{center}
\end{figure}

The graph differential makes the graphs into a chain complex $\Gamma_{\bullet}$, at degree $n$ the space
  $\Gamma_n$ consists of linear combinations of graphs with $n$ vertices. The graph homology is defined as the quotient of graph cycles by graph boundaries in the usual manner
\bea H_n=\frac{Z_n}{B_n}~.\nn\eea
And dually, we can define the graph cochain complex $\Gamma_n^*$ whose $n$-cochains are linear mappings $c^n:\,\Gamma_n\to \BB{R}$. We usually write this linear map as a pairing
\bea c^n\circ c_n=\bra c^n,c_n\ket~;~~~c^n\in\Gamma^*_n~,~~c_n\in \Gamma_n~.\nn\eea
We present the Figure \ref{ex_box_diff}
as an illustration of the graph differential. Here we use the oriented labelled graphs as representative of equivalence class.
 Since under the equivalence relation the edges starting and ending at the same vertex are forbidden then it is clear that
  we get the factor $6$ for the first line  and the factor $2$ for the second line.  Later on we will come back to the algebraic description
   of the graph differential coming from the perturbative treatment of integrals.

\subsection{Kontsevich Theorem}
\label{KT}

In this subsection we review briefly the Kontsevich theorem \cite{Kontsevich:symplectic,  Kontsevich:Feynamn}.
 Here our presentation is quite formal and we present
 the theorem as certain ad hoc recipe. Later in section \ref{bv-isomorphism} we will rederive this theorem in much more natural fashion.

The Kontsevich theorem is about the isomorphism between the graph complex and certain Chevalley-Eilenberg complex.
The Chevalley-Eilenberg complex has been reviewed in section \ref{AVOTI}.
Here we consider the Lie algebra of formal Hamiltonian vector fields on $\BB{R}^{2n}$ equipped with the constant symplectic
 structure $\Omega_{\mu\nu}$. The generalization to $\BB{R}^{2n|m}$ is straightforward and we leave it aside for the moment.
  These vector fields are generated by formal polynomial functions on $\BB{R}^{2n}$
\bea \BB{X}_f=(\partial_{\mu}f(x))(\Omega^{-1})^{\mu\nu}\partial_{\nu}~.\nn\eea
The formal Hamiltonian vector field form a closed algebra due to the relation $[\BB{X}_f,\BB{X}_g]=\BB{X}_{\{f,g\}}$.
 We are interested in the Lie algebra $\textrm{Ham}^0 (\BB{R}^{2n}) $ which contains  polynomials that have only quadratic and higher terms.
 In particular  we are interested in  the Chevalley-Eilenberg complex for  $\textrm{Ham}^0 (\BB{R}^{2n}) $.
  Let us make the following
 abbreviation for the chains in this CE complex
\bea
\BB{X}_{f_1}\wedge... \wedge\BB{X}_{f_k}\Rightarrow (f_1,...,  f_k)\nn
\eea
 with the boundary operator defined as usual
\bea
\partial(f_1,...,  f_k)=\sum_{i<j}(-1)^{i+j+1}(\{f_i,f_j\},f_1, ... , \hat{f_i},... , \hat{f_j}, ... , f_k)~.\label{def-defCE-comp}\eea
    This defines the  Chevalley-Eilenberg complex for  $\textrm{Ham}^0 (\BB{R}^{2n}) $ which we denote ${\rm CE}_{\bullet}(\textrm{Ham}^0(\BB{R}^{2n}))$.
     Obviously  we define the dual objects, the cochains evaluating on a chain as
\bea c^k (f_1,\cdots f_k)~\in ~\BB{R}~.\nn\eea

 The Kontsevich's theorem\footnote{In fact, Kontsevich studied a relative complex ${\rm CE}_{\bullet}(\textrm{Ham}^0(\BB{R}^{2n});\FR{sp}(2n))$ which would correspond to graphs with cubic and higher vertices and it has a larger cohomology group. But as far as the mapping (\ref{iso_CE_gph1})
   is concerned there is no difference.} states  that there is an isomorphism between two complexes
\bea
\Gamma_{\bullet}\sim \lim_{n\to \infty}{\rm CE}_{\bullet}(\textrm{Ham}^0(\BB{R}^{2n}))~,\label{iso_CE_gph1}\eea
where  $\Gamma_{\bullet}$ is graph complex defined in subsection \ref{BOGT}.
For any fixed $n$, the mapping is a homomorphism, meaning it maps the graph differential into the  Chevalley-Eilenberg
  differential. Taking the limit $n\to\infty$ is for the sake of accommodating graphs with arbitrary number of vertices.

We will now give the formal recipe of the mapping (\ref{iso_CE_gph1}) in detail.  Given a chain $(f_1, ... ,f_k)$ from
 ${\rm CE}_{\bullet}(\textrm{Ham}^0(\BB{R}^{2n}))$ we can construct the graph chain from $\Gamma_\bullet$
   following steps:
\begin{itemize}
 \item Fix vertices numbered from 1 to $k$, the vertex $i$ can be maximally $p_i$-valent if the polynomial
  $f_i$ has highest degree $p_i$,
  all vertices are considered distinct but legs from one vertex are considered identical.
 \item Exhaust all possible unoriented graphs by drawing edges between vertices, one gets a collection of graphs. There are only finite number of graphs since each vertex has a maximum valency and at this stage we allow graphs containing 1- or 2-valent vertices.
\begin{itemize} \item For each graph one assigns an arbitrary orientation to its edges, then for an edge from vertex $i$ to $j$, one differentiates $f_i$ with $\partial_{\mu}$, $f_j$ with $\partial_{\nu}$ and contracts $\mu,\nu$ with $(\Omega^{-1})^{\mu\nu}$.
 \item Do this for all edges, then set all $x$'s to zero, the result is a c-number.
 \item Divide the c-number by $\#P$ (defined in subsection \ref{IiRnGIaFD}), this is the coefficient of the graph.
 \end{itemize}
 \item Doing this for all graphs one gets a linear combination of graphs.
 \item Combine coefficients of the graphs that belong to the same equivalence class (see subsection \ref{BOGT})
  and this give us the graph chain
\end{itemize}
We remark that step 2 can be greatly simplified, but the current low pace approach has a closer link to perturbation theory. We will give the simplified version after the next example.

The hard part is to show that the recipe gives a homomorphism, i.e. the differential in ${\rm CE}_{\bullet}(\textrm{Ham}^0(\BB{R}^{2n}))$
 is mapped to the graph differential in $\Gamma_\bullet$ under this recipe.  This can be proved using methods from \cite{ConantVogtmann}.
However there is  striking similarity between the recipe given above and the Feynman diagrams in perturbative treatment of integrals.
 This may suggest that a simpler proof can arise from a carefully designed Gaussian integral, whose perturbative expansion naturally gives the above recipe. But there is a technical difficulty: we need to connect legs of the vertices using the inverse symplectic structure $\Omega^{-1}$,
but we cannot possibly write a Gaussian integral with  $\exp\{-x^{\mu}\Omega_{\mu\nu}x^{\nu}\}$ since it is trivial.
 But let us put aside this problem for now and give instead the concrete example of the recipe where we will hopefully clarify some of
  the steps.

\subsubsection{Example}\label{example-KT}

\begin{figure}[h]
\begin{center}
\psfrag{G1}{\scriptsize{$\Gamma_1$}}\psfrag{G2}{\scriptsize{$\Gamma_2$}}
\psfrag{G3}{\scriptsize{$\Gamma_3$}}\psfrag{G4}{\scriptsize{$\Gamma_4$}}\psfrag{G5}{\scriptsize{$\Gamma_5$}}
\psfrag{G6}{\scriptsize{$\Gamma_6$}}\psfrag{G7}{\scriptsize{$\Gamma_7$}}
\includegraphics[width=3.8in]{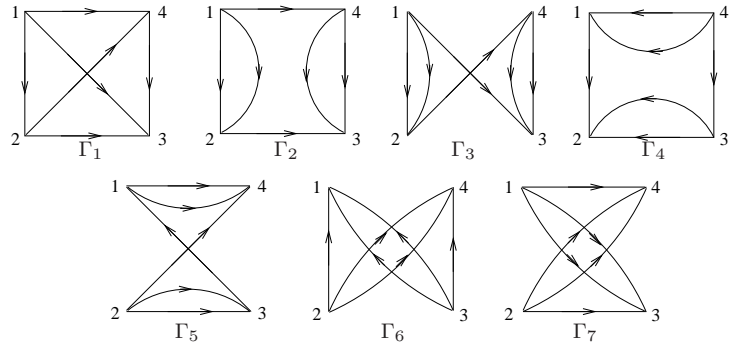}
\caption{Wick's contraction}\label{ex_box_fig}
\end{center}
\end{figure}

Let $e,f,g,h$ be four cubic functions on $\BB{R}^{2n}$ with the normalization given by (\ref{normalization-pol}).
 Let us first fix $e,f,g,h$ to be on the vertices numbered 1,2,3,4.  There are seven different ways of
connecting the legs together as in Figure \ref{ex_box_fig}.
In fact, the last six oriented labelled graphs  are related by relabeling of vertices.
Following the recipe we can construct  for each graph the following  coefficient $c_{\Gamma_i}$
\bea &&c_{\Gamma_1}=X(e,f,g,h)=e^{\ga\lambda\gd}f_{\ga}^{~\gb\kappa}g_{\gb\gc\lambda}h^{\gc}_{~\gd\kappa}~,
~~~c_{\Gamma_2}=\frac14Y(e,f,g,h)=\frac14e^{\ga\lambda\gd}f_{\ga~\lambda}^{~\gb}g_{\gb\gc\kappa}h^{\gc~\kappa}_{~\gd}~,\nn\\
&&c_{\Gamma_3}=\frac14Y(e,f,h,g)~,~~c_{\Gamma_4}=\frac14Y(h,e,f,g)~,~~c_{\Gamma_5}=\frac14Y(f,g,e,h)~,\nn\\
&&c_{\Gamma_6}=-\frac14Y(f,h,g,e)~,~~c_{\Gamma_7}=-\frac14Y(e,g,f,h)~.\nn\eea
In these expressions, instead of writing $\Omega^{-1}$ explicitly, we have raised indices using the inverse of symplectic structure, i.e.
 $(...)^{\mu}=(...)_{\rho}(\Omega^{-1})^{\rho\mu}$. And we also denote $f_{\mu\nu\rho}=\partial_{\mu}\partial_{\nu}\partial_{\rho}f$.
 Thus finally, the graph chain is
\bea \Gamma&=&\Gamma_1\cdot X(e,f,g,h)+\Gamma_2\cdot\frac14\Big(Y(e,f,g,h)-Y(e,f,h,g)\nn\\
&&\hspace{3cm}-Y(h,e,f,g)+Y(f,g,e,h)+Y(f,h,g,e)-Y(e,g,f,h)\Big)\nn\\
&=&\frac1{24}\Gamma_1\Big(X(e,f,g,h)+\textrm{asym perm's}\Big)
+\frac14\Gamma_2\cdot\Big(\frac14Y(e,f,g,h)+\textrm{asym perm's}\Big)~.\label{gph_ch}\eea
The two factors 1/24 and $1/4\cdot1/4$ are recognized as the total symmetry factor $\#P\#V$ of $\Gamma_{1,2}$ (in this case loops are clearly forbidden). This is the graph chain associated with the  ${\rm CE}_{\bullet}(\textrm{Ham}^0(\BB{R}^{2n}))$  chain $(e,f,g,h)$.

We also see that there is a slight shortcut to the recipe above, we need only do the first two sub-steps in step 2 for one \emph{representative} of an equivalence class of graphs and manually sum over anti-symmetric permutations of $f_i$, then divide the result by the total symmetry factor of the graph.

Next we apply the graph differential to the chain $\Gamma$ defined in (\ref{gph_ch}).
 We denote the bottom right graph in Figure \ref{ex_box_diff} by $\Gamma_8$.
Then graph differential applied to (\ref{gph_ch}) gives us
\bea
\partial \Gamma =\Gamma_8(-\frac14X(e,f,g,h)+\frac12\frac14Y(e,f,g,h)+\textrm{asym-perms})~.\label{gph_ch_diff}\eea

Here comes the key moment, we need to compare this result to the result obtained by applying  the Chevalley-Eilenberg
 differential $\partial_{CE}$
  to $(e,f,g,h)$ first and the recipe above second.
To do this, we write out $X$ and $Y$ in (\ref{gph_ch_diff})
\bea \partial \Gamma
&=&\Gamma_8\Big(f^{\ga\gb\kappa}g_{\gb}^{~\gc\lambda}\big(\frac14e_{\ga\lambda}^{~~\gd}h_{\gd\gc\kappa}
+\frac18e_{\ga\kappa}^{~~\gd}h_{\gd\gc\lambda}\big)+\textrm{asym-perms}\Big)\nn\\
&=&\Gamma_8\Big(f^{\ga\gb\kappa}\frac12g_{\gb}^{~\gc\lambda}\big(\frac12e_{\ga\lambda}^{~~\gd}h_{\gd\gc\kappa}
+\frac18e_{\ga\kappa}^{~~\gd}h_{\gd\gc\lambda}+\frac18h_{\ga\kappa}^{~~\gd}e_{\gd\gc\lambda}\big)+\textrm{asym-perms}\Big)\nn\\
&=&\Gamma_8\Big(\frac{1}{16}f^{\ga\gb\kappa}g_{\gb}^{~\gc\lambda}\{e,h\}_{\ga\kappa\gc\lambda}+\textrm{asym-perms}\Big)~.\label{gph_ch_diff_full}\eea
 The Chevalley-Eilenberg
 differential $\partial_{CE}$   acting on $(e,f,g,h)$ gives
\bea \partial_{CE}(e,f,g,h)&=&(\{e,f\},g,h)-(\{e,g\},f,h)+ ... + (\{g,h\},e,f)\nn\\
&=&\frac14\big((\{e,f\},g,h)+\textrm{asym-perms}\big)~.\nn\eea
Going back to recipe the first term gives rise to the following graph chain
\bea (\{e,f\},g,h)\to\Gamma_8\frac14\{e,f\}^{\ga\gb\gc\gd}g_{\gc\gd}^{~~\,\kappa}h_{\kappa\ga\gb}\nn\eea
 and so on.  Finally we get the following correspondence
\bea \partial_{CE}(e,f,g,h)\to \Gamma_8\frac1{16}\Big(\{e,f\}^{\ga\gb\gc\gd}g_{\gc\gd}^{~~\,\kappa}h_{\kappa\ga\gb}
+\textrm{asym-perms}\Big)~,\nn\eea
which is in full agreement with (\ref{gph_ch_diff_full}) where we applied the graph differential explicitly.

\subsection{Algebraic Description of Graph Chains}\label{ADoGC}

Due to the combinatorial nature of the graph complex, the proof of any proposition involving graph complex is extremely cumbersome.
 Therefore  here  we suggest a more algebraic description of the graph complex, which is motivated by the perturbative expansion
  of integrals.
 The construction involves introducing some formal parameters to represent edges and vertices in such a way that a graph corresponds to a polynomial in these parameters. These parameters have to conform to the symmetry properties of a graph, which prompts us to use odd variables $t_i,i=1\cdots n$ ($n$ is the number of vertices) to represent the vertices, and even variables $t_{ij}$  to represent edges.
  We assume that $t_{ij}=-t_{ji}$ in order to take the orientation of the edge into account.

  Every oriented labelled graph with $N$ vertices can be represent by a monomial according to the following prescription

   \begin{itemize}
 \item  for every vertex we include $t_i$, so from all vertices we have $t_1 t_2 ... t_n$;
 \item for each oriented edge starting from vertex $i$  and ending at vertex $j$, include a factor $t_{ij}$.
 \end{itemize}

 Let us illustrate these rules by some examples. If we look at $\Gamma_1$ from Figure \ref{ex_box_fig}
  as concrete oriented labelled graph  then there is the following monomial
  \bea
   t_1 t_2 t_3 t_4~ t_{14} t_{13} t_{12} t_{24} t_{23} t_{43}~.\nn
  \eea
 If we multiply this monomial by $-1$ then we will change the orientation of the graph.  Another example is
  the graph $\Gamma_2$ from Figure \ref{ex_box_fig}. If $\Gamma_2$ is understood as oriented labelled graph
   then there is the following monomial
\bea
 t_1 t_2 t_3 t_4~ t_{12}^2 t_{43}^2 t_{14} t_{23}~.\nn
\eea

   However we are interested in the equivalence classes of graphs as described in subsection \ref{BOGT}.
    The equivalence classes can be represented by polynomials where we sum over all equivalent oriented labelled graphs
     (monomials).
   To write down a polynomial representing a graph, we follow the steps
 \begin{itemize}
 \item assign indices $l_1, ... , l_n$ to all  $n$ vertices of a graph;
 \item for each edge from vertex $i$ to $j$, include a factor $t_{l_il_j}$;
 \item multiply the polynomial by $t_{l_1}\cdots t_{l_n}$ and sum $l_1\cdots l_n$ from 1 to $n$.
 \end{itemize}
For example, the graph  $\Gamma_1$ in Figure \ref{ex_box_fig} understood as representative of the equivalence class of graphs
gives the following polynomial
\bea
\sum\limits_{l_1=1}^4 \sum\limits_{l_2=1}^4\sum\limits_{l_3=1}^4 \sum\limits_{l_4=1}^4 ~t_{l_1}t_{l_2}t_{l_3}t_{l_4}~
t_{l_1l_2}t_{l_1l_4}t_{l_1l_3}t_{l_2l_4}t_{l_4l_3}t_{l_2l_3}~.\nn\eea
 While the graph $\Gamma_2$ in Figure \ref{ex_box_fig} understood as representative of the equivalence class of graphs
  gives
  \bea
 \sum\limits_{l_1=1}^4 \sum\limits_{l_2=1}^4\sum\limits_{l_3=1}^4 \sum\limits_{l_4=1}^4~  t_{l_1} t_{l_2} t_{l_3}
  t_{l_4}~ t_{l_1l_2}^2 t_{l_4l_3}^2 t_{l_1 l_4} t_{l_2 l_3}~.\nn
  \eea
Notice that these polynomials \emph{do not contain any symmetry factors}.
The point of summing over the dummy indices $l_1\cdots l_n$ is so that the symmetry group $s_n$ acts trivially on the polynomial. Hence one such polynomial represents an \emph{equivalence class of graphs}. In this formalism the orientation is taken into account automatically.
 Thus we can think about the graph chains as the formal polynomials in odd parameters $t_i$ and even parameters $t_{ij}$.
  Later on we will present the graph differential as formal differential operator acting on these polynomials. Many calculations
   drastically simplify with this algebraic description of graphs.

Sometime we will need to solve the following problem. Given polynomial in $t_i$'s and $t_{ij}$'s  we would like to know
 if a concrete equivalence class of graphs is present in this polynomial and if yes then which numerical factor is in front of this graph.
While given a polynomial, to recover concrete equivalence class of graphs it represents we just need to do the following
\begin{itemize}
\item Pick an arbitrary graph out of an equivalence class of graphs, say, a representative whose vertices are labeled $1, ...,  n$;
\item differentiating with respect to the odd parameters $\partial_{t_n}\cdots \partial_{t_1}$;
\item for each edge from vertex $i$ to $j$, include a derivative $\partial_{t_{ij}}$;
\item set to zero all formal parameters;
\item divide by the symmetry factor $\#V\#P$.
\end{itemize}
Observe that no summation over the dummy indices is needed. After this procedure we end up with numerical coefficient
 which will tell us if the equivalence class is present or not and with what factor.  For example, to see if the equivalence class
  of $\Gamma_1$ from Figure \ref{ex_box_fig}  is present we have to apply the following operator to the polynomial
\bea  \frac{1}{24}\partial_{t_4}\cdots \partial_{t_1}\partial_{t_{12}}\partial_{t_{14}}
\partial_{t_{13}}\partial_{t_{24}}\partial_{t_{43}}\partial_{t_{23}} \nn
\eea
 and then set all parameters to zero.  For the equivalence class represented by the graph $\Gamma_2$ from Figure \ref{ex_box_fig}
  the corresponding differential operator is
\bea
~~\frac{1}{4}\frac{1}{4}\partial_{t_4}\cdots \partial_{t_1}\partial^2_{t_{32}}\partial^2_{t_{41}}
\partial_{t_{12}}\partial_{t_{43}}~.\nn\eea
Note  the symmetry factors for $\Gamma_1$ is $\#V=24,\#P=1$,
 and for $\Gamma_2$ are $\#V=4,\#P=4$.

We will derive later the operator corresponding to the graph differential and discuss more the algebraic description of graphs
 in the context of perturbative expansion.

\section{BV formalism and graph complex}
 \label{bv-isomorphism}

 This section presents the non-trivial application of finite dimensional BV formalism.   We will
  reprove the  Kontsevich theorem about the relation between graph complex $\Gamma_\bullet$ and Chevalley-Eilenberg
   complex ${\rm CE}_{\bullet}(\textrm{Ham}^0(\BB{R}^{2n}))$. The existing proofs can be found in \cite{ConantVogtmann,Hamilton-2005} and
    in the appendix of \cite{2010arXiv1006.1240Q}  where we gave the  proof generalizing Kontsevich's theorem to the case of chord diagrams and extended Chevalley-Eilenberg complexes suited for the study of knots.  The proof presented here is streamlined
     and simplified version of the proof from  \cite{2010arXiv1006.1240Q}.

\subsection{A Universal BV Theory on a Lattice}
\label{AUBVToaL}

We mentioned in subsection \ref{KT} that the mapping from   the Chevalley-Eilenberg complex to graph complex is nothing but the application of Feynman rules and that there is a technical difficulty in realizing the naive Feynman rules. We would like to construct a BV theory that circumvents this difficulty and furthermore, whose path integral gives, instead of numbers, graphs as outcome.
 In this way, it turns out that the Ward identity directly imply  Kontsevich's theorem. We will embed into BV formalism
  the perturbation theory presented in subsection \ref{IiRN}.

 Let us assume that the vector space $\BB{R}^{2n}$ is equipped with constant symplectic structure $\Omega_{\mu\nu}$
  and as in  subsection \ref{IiRN}  we consider $N$ copies of this vector space, $\bigoplus\limits_{i=1}^N \mathbb{R}^{2n}$.
   We use $x^\mu_i$ to denote the coordinates in this big vector space, where $\mu=1,2,..., 2n$ and $i=1,2,....N$.
     Now let us construct new graded manifold
     $$\bigoplus\limits_{i=1}^N \mathbb{R}^{2n} \oplus  \bigoplus\limits_{i=1}^N \mathbb{R}^{2n}[-1]~,$$
    where in addition to $x^\mu_i$ we introduced the odd coordinates $\xi^\mu_i$ of degree $-1$.  On the space of
     functions $C^\infty (\bigoplus\limits_{i=1}^N \mathbb{R}^{2n} \oplus  \bigoplus\limits_{i=1}^N \mathbb{R}^{2n}[-1])$ we
      can introduce the structure of BV-algebra (see the definitions \ref{BValgebra}  and \ref{BValgebra-alternative}).
 As in section \ref{ss-odd-FT} the odd Laplacian operator is defined as follows
 \bea
 \Delta=\sum_{i=1}^N~(\Omega^{-1})^{\mu\nu}\frac{\partial}{\partial \xi_i^{\mu}}\frac{\partial}{\partial x_i^{\nu}}~.\nn\eea
 The odd Poisson bracket is defined accordingly as
 \bea
 \{ g, h \} =\sum_{i=1}^N (\Omega^{-1})^{\mu\nu}  \left ( \frac{\partial  g}{\partial x_i^\mu} \frac{\partial h}{\partial \xi^\nu_i} + (-1)^{|g|} \frac{\partial g}{\partial \xi^\nu_i} \frac{\partial h}{\partial x^\mu_i} \right )~,
 \eea
  where $g, h \in C^\infty (\bigoplus\limits_{i=1}^N \mathbb{R}^{2n} \oplus  \bigoplus\limits_{i=1}^N \mathbb{R}^{2n}[-1])$.
 If we introduce $\xi_{i\mu} = \Omega_{\mu\nu}\xi^\nu_i$ then we deal with the odd cotangent bundle $T^*[-1] M$, where $M= \bigoplus\limits_{i=1}^N \mathbb{R}^{2n}$.  Therefore all our previous discussion from section \ref{oddFT} is applicable here.

 Now let us construct some specific functions on  $C^\infty (\bigoplus\limits_{i=1}^N \mathbb{R}^{2n} \oplus  \bigoplus\limits_{i=1}^N \mathbb{R}^{2n}[-1])$.
Let us pick a function $f(x)$ on $\BB R^{2n}$ and auxiliary odd parameters
  $t_i$  of degree $-1$ ($i=1,2, ..., N$). We can construct the following function on our BV manifold
\bea
{\cal O}[f] =\sum_{i=1}^N~ \big(t_if(x_i)+\xi_i^{\mu}\partial_{\mu}f(x_i)\big)~.\nn\eea
 It is quite easy to see that $\Delta {\cal O} [f]=0$, therefore the bracket between ${\cal O}[f]$ and ${\cal O}[g]$ is induced as
\bea
\{ {\cal O} [f], {\cal O}[g]\}=-\Delta({\cal O}[f] {\cal O}[g])~.\nn\eea
We can directly calculate the bracket $\{{\cal O}[f], {\cal O}[g]\}$
\bea
\{{\cal O}[f], {\cal O}[g]\}&=&\sum_i-\big(\partial_{\rho}f(x_i)\big)\Omega^{\sigma\rho}\partial_{\sigma}\big(t_ig(x_i)+\xi_i^{\gc}\partial_{\gc}g(x_i)\big)\nn\\
&&\hspace{3cm}+
\sum_i\partial_{\sigma}\big(t_if(x_i)+\xi_i^{\gc}\partial_{\gc}f(x_i)\big)\Omega^{\sigma\rho}\big(\partial_{\rho}g(x_i)\big)\nn\\
&=&\sum_i~\Big(2t_i\{f,g\}(x_i)+\xi^{\sigma}_{i}\partial_{\sigma}\{f,g\}(x_i)\Big)~,\label{unequal}\eea
 where the bracket $\{f, g\}$ is the standard Poisson bracket on $\BB{R}^{2n}$ with respect to $\Omega$. Let us point out  that
  $\{{\cal O}[f], {\cal O}[g]\} \neq 2 {\cal O}[\{ f, g\}]$.  Next we introduce the special function which
   is the one we used in (\ref{Gauss_feyn_high_d1}),
\bea
S = \frac12\sum\limits_{i,j=1}^N  x_i^\mu t^{ij} x_j^\nu ~\Omega_{\mu\nu}~,\label{BV_action_discrete}\eea
 which physicists call BV action. Here
$t_{ij},~i,j=1,\cdots N$ are now formal degree 0 parameters with $t_{ij}=-t_{ji}$ and $t^{ik}t_{kj}=\delta^i_j$.
 The purpose of introducing these parameters is to make our perturbation theory universal in the sense that the
   Feynman diagrams are computed as a function of $t_{ij}$.

   Now we are ready to state and reprove the Kontsevich theorem discussed in subsection \ref{KT}.  Below we show that
    the theorem is simple consequence of BV formalims.

   \begin{Thm}[Kontsevich]\label{theorem-cocycle}
    Let us take a collection of polynomial functions $f_1, ... f_l$ on $\BB{R}^{2n}$ and define the following integral
\bea \bra (f_1, ... , f_l) \ket=\frac{1}{Z[0]}\int d^{2n}x_1\cdots d^{2n}x_N
~\sum_{i_1=1}^N t_{i_1}f_1(x_{i_1}) ... \sum_{i_l=1}^N t_{i_l}f_l(x_{i_l})~e^{-S}~\in \BB{R}[t_i,t_{ij}]~,\label{inegarakkee22}\eea
where as answer we obtain a polynomial in $t_i$ and $t_{ij}$ which can be understood as a graph chain
  described in subsection \ref{ADoGC}.  This integral
 can be understood as map from chain $(f_1, ... , f_l)$ in ${\rm CE}_{\bullet}(\textrm{Ham}^0(\BB{R}^{2n}))$ to graph chain
  in $\Gamma_\bullet$.
 The reader can compare the perturbative expansion of the formula (\ref{inegarakkee22})  to the recipe in section \ref{KT}
  and indeed they are the same. Moreover the map (\ref{inegarakkee22}) satisfies the property
\bea
\bra\partial (f_1,\cdots f_l)\ket=\partial_{Gph}\bra(f_1,\cdots f_l)\ket~,\label{iso_CE_gph-1}
\eea
where $\partial$ is the Chevalley-Eilenberg differential define in  (\ref{def-defCE-comp})   and $\partial_{Gph}$ is the graph differential,
 which has the following explicit form\footnote{Here the cumbersome factor $(N-l+1)$ arises because we chose the number of lattice cites $N$ to be a large number and therefore we can accommodate graphs of varying number of vertices at the same time. In fact, for the current problem, one may well take $N=l$ and thereby eliminate this ugly factor.}
\bea
\partial_{Gph}=-\frac{1}{2(N-l+1)}\sum_{k,p=1;k\neq p}^NR^{p}_{k}\frac{\partial}{\partial t_{kp}}\frac{\partial}{\partial t_{p}}~,\label{gph_diff_poly}\eea
 when we deal with graph chains as polynomials (the precise definition of the operator $R^{p}_{k}$ see below in the proof).
\end{Thm}

Now we provide the proof and further explanation for this theorem.  We leave to the reader to check that
 the perturbative expansion of (\ref{inegarakkee22}) with the propagator
\bea
\bra x_i^{\mu} x_j^{\nu}\ket=(\Omega^{-1})^{\mu\nu}t_{ij} \nn\eea
coincides with the prescription given in subsection \ref{KT} when we understand the answer as graph chain (i.e., formal
 polynomial in  $t_i$'s and $t_{ij}$'s).

 We concentrate on the proof of the relation (\ref{iso_CE_gph-1}).  For this we have to embed the integral (\ref{inegarakkee22}) into
  the BV framework.  Let us introduce short hand notation for the integration measure ${\cal D} x = d^{2n}x_1\cdots d^{2n}x_N$.
   The integral (\ref{inegarakkee22}) can be thought as integral over odd conormal bundle
\bea \bra (f_1, ..., f_l) \ket=\int\limits_{\xi^=0}{\cal D}x
~{\cal O} [f_1] {\cal O}[f_2] ... {\cal O}[f_l]~e^{-S}~,\nn\eea
 where the functions ${\cal O}[f]$ and $S$ are defined above on the BV manifold $\bigoplus\limits_{i=1}^N \mathbb{R}^{2n} \oplus  \bigoplus\limits_{i=1}^N \mathbb{R}^{2n}[-1]$. Using the terminology from subsection \ref{s-integration} here the submanifold $C$ coincides with $M$
  and thus there is no odd conormal directions left.
As explained earlier in subsection \ref{s-integration}, there is one single Ward identity in BV formalism (\ref{Ward1}): $\int\Delta(\cdots)=0$, we analyze its implication in the present setting
\bea
\int\limits_{\xi=0}{\cal D}x
~\Delta\big({\cal O}[f_1] {\cal O}[f_2] ...  {\cal O}[f_l]~e^{-S}\big)=0~.\nn\eea
The standard manipulation with the odd  Laplacian operator (\ref{Delta-defPB}) leads to
\bea &&\Delta\big({\cal O}[f_1]{\cal O}[f_2] ...  {\cal O} [f_l]~e^{-S}\big)\nn\\
&=&\Delta\big({\cal O}[f_1]{\cal O}[f_2] ... {\cal O}[f_l]\big)~e^{-S}+(-1)^ l {\cal O}[f_1] {\cal O}[f_2] ...  {\cal O}[f_l]\Delta e^{-S}\nn\\
&&-(-1)^l\big\{{\cal O}[f_1]{\cal O}[f_2] ... {\cal O}[f_l],S\big\}e^{-S}~.\label{interm1}\eea
We observe that $\Delta S=0$ trivially. And furthermore all the above expressions will be restricted onto the submanifold given by $\xi=0$ and thus we can make the replacement
\bea\Delta\big({\cal O}[f_1]{\cal O}[f_2] ... {\cal O}[f_l]\big)\Big|_{\xi=0}=-2\sum_{i<j}(-1)^{i+j+1}{\cal O}[\{f_i,f_j\}] {\cal O}[f_1] ... \widehat{{\cal O}[f_i]}... \widehat{{\cal O}[f_j]} ...  {\cal O}[f_l]\Big|_{\xi=0}~,\nn\eea
 where we have used the relation (\ref{unequal}).
Thus we conclude that the first term of  (\ref{interm1}) gives the correlator
\bea \sum_{i<j}(-1)^{i+j+1}\bra(\{f_i,f_j\},f_1, ..., \widehat{f_i}, ...,\widehat{f_j},..., f_l )\ket=\bra \partial(f_1,..., f_l)\ket~,\nn\eea
i.e. the correlator of the Chevalley-Eilenberg differential of $(f_1,\cdots f_n)$.
 Next we analyze the last term of  (\ref{interm1})
\bea &&-(-1)^l\int\limits_{\xi=0}{\cal D}x~\big\{{\cal O}[f_1]{\cal O}[f_2] ...  {\cal O}[f_l],S\big\}e^{-S}\nn\\
&=&-\sum_{p=1}^l (-1)^p \int\limits_{\xi=0}{\cal D}x~{\cal O}[f_1] ... \{{\cal O}[f_p],S\} ...  {\cal O}[f_l]~e^{-S}\nn\\
&=&\sum_{p=1}^l (-1)^p \int\limits_{\xi=0}{\cal D}x~{\cal O}[f_1] ... \Big(\underbrace{-\sum_{i_p,j}~t^{i_pj}\partial_{\mu}f_p(x_{i_p})x^{\mu}_j}_{p^{th}}\Big)... {\cal O}[f_l]~e^{-S}~.\label{interm}\eea
Now it is best that we should explain where we are heading before we make the plunge. Focusing on the round brace, assume that $x_j$ is connected to $f_k$ in the ensuing Gaussian integral, we will get a factor
\bea \sum_j\sum_{i_k}~(\Omega^{\mu\nu}t_{ji_k})(t_{i_k}\partial_{\nu}f_k(x_{i_k}))\Big(-\sum_{i_p}~t^{i_pj}\partial_{\mu}f_p(x_{i_p})\Big)~,\nn\eea
where the quantity in the first brace comes from the propagator. Summing over $j$, we get $\delta^{i_p}_{i_k}$, which renames all the $i_p$ into $i_k$. This leads of course to
\bea \sum_{i_k}~(\Omega^{\mu\nu})~t_{i_k}\partial_{\nu}f_k(x_{i_k})(\partial_{\mu}f_p(x_{i_k}))
=\sum_{i_k}~t_{i_k}\{f_p,f_k\}(x_{i_k})~.\nn\eea

The above describes exactly the process of combining two vertices by shrinking an edge between them and letting the new vertex inherit all the other edges belonging to the two old ones. However it will be more convenient if we could write the effect of $\partial_{Gph}$ as an operator acting on the polynomial $\BB{R}[t_{ij},t_i]$. Thus instead we proceed from (\ref{interm}) and perform the Gaussian integral by replacing $x_i^{\mu}$ with $\partial_{J^i_{\mu}}$ as in (\ref{Gauss_feyn_high_d})
%
%
\bea (\ref{interm})&=&\sum_{i_1\cdots i_l=1}^N\sum_{p=1}^l (-1)^p \big(t_{i_1}f_1(\frac{\partial}{\partial J^{i_1}})\big)\cdots \nn\\
&&\hspace{3cm}\cdots \Big(\underbrace{-\sum_{j=1}^N~t^{i_pj}\partial_{\mu}f_p(\frac{\partial}{\partial J^{i_p}})\frac{\partial}{\partial J_{\mu}^j}}_{p^{th}}\Big)\cdots \big(t_{i_l}f_l(\frac{\partial}{\partial J^{i_l}})\big)~e^{\frac12J\Omega^{-1}tJ}\Big|_{J=0}~,\nn\eea
 where $J\Omega^{-1}tJ$ is the short hand notation for $\sum\limits_{i,j=1}^NJ^i_\mu (\Omega^{-1})^{\mu\nu} t_{ij} J^j_\nu$.
Letting $\partial_{J_{\mu}^j}$ hitting the last exponential
\bea
&&\sum_{i_1\cdots i_l=1}^N\sum_{p=1}^l (-1)^p \big(t_{i_1}f_1(\frac{\partial}{\partial J^{i_1}})\big)\cdots\nn\\
&&\hspace{1cm}\cdots \Big(\underbrace{-\sum_{j=1}^N~t^{i_pj}\partial_{\mu}f_p(\frac{\partial}{\partial J^{i_p}})}_{p^{th}}\Big)
\cdots \big(t_{i_l}f_l(\frac{\partial}{\partial J^{i_l}})\big)
~\sum_{k=1}^N(\Omega^{-1})^{\mu\nu}t_{jk}J^k_{\nu}~e^{\frac12J\Omega^{-1}tJ}\Big|_{J=0}\nn\\
&=&\sum_{i_1\cdots i_l=1}^N\sum_{p=1}^l (-1)^p \big(t_{i_1}f_1(\frac{\partial}{\partial J^{i_1}})\big)\cdots\nn\\
&&\hspace{1cm}\cdots \big(\underbrace{-\partial_{\mu}f_p(\frac{\partial}{\partial J^{i_p}})}_{p^{th}}\big)
\cdots \Big(t_{i_l}f_l(\frac{\partial}{\partial J^{i_l}})\Big)~(\Omega^{-1})^{\mu\nu}J^{i_p}_{\nu}~e^{\frac12J\Omega^{-1}tJ}\Big|_{J=0}~.\label{star}\eea

We claim that this expression can be written as
\bea -\frac{1}{N-l+1}\sum_{k,q=1}^NR^q_k\frac{\partial}{\partial t_{kq}}\frac{\partial}{\partial t_{q}}\int_{\xi=0}{\cal D}x~{\cal O}_1\cdots {\cal O}_l~e^{-S}~,\label{stella}\eea
where $R^p_q$ is a renaming operator acting on the polynomials of $t_{ij}$ that renames $p$ to $q$
\bea R^p_q t_{ij}=\bigg\{\begin{array}{cc}
                           t_{ij} & i,j\neq p \\
                           t_{qj} & i=p \end{array}~.\label{rename}\eea
To see this, again we replace $x$ with $\partial_{J}$
\bea (\ref{stella})&=&-\sum_{k,q=1}^NR^q_k\frac{\partial}{\partial t_{kq}}\frac{\partial}{\partial t_{q}}\bigg(\sum_{i_1\cdots i_l=1}^N\big(t_{i_1}f_1(\frac{\partial}{\partial J^{i_1}})\big)\cdots \big(t_{i_l}f_l(\frac{\partial}{\partial J^{i_l}})\big)~e^{\frac12J\Omega^{-1}tJ}\Big|_{J=0}\bigg)\nn\\
&=&-\sum_{k,q=1}^NR^q_k\frac{\partial}{\partial t_{q}}\bigg(\sum_{i_1\cdots i_l=1}^N\big(t_{i_1}f_1(\frac{\partial}{\partial J^{i_1}})\big)\cdots \big(t_{i_l}f_l(\frac{\partial}{\partial J^{i_l}})\big)~(J^k_{\rho}(\Omega^{-1})^{\rho\sigma}J_{\sigma}^q)e^{\frac12J\Omega^{-1}tJ}\Big|_{J=0}\bigg)~.\nn\eea
We now commute $J_{\sigma}^q$ to the left most position, picking up only a commutator
\bea
&&-\sum_{k=1}^N\sum_{p=1}^l\sum_{i_1\cdots i_l=1}^N R^{i_p}_k\frac{\partial}{\partial t_{i_p}}\bigg(\big(t_{i_1}f_1(\frac{\partial}{\partial J^{i_1}})\big)\cdots
\Big(t_{i_p}\partial_{\sigma}f_p(\frac{\partial}{\partial J^{i_p}})\Big)\cdots\nn\\
&&\hspace{4cm}\cdots
\big(t_{i_l}f_l(\frac{\partial}{\partial J^{i_l}})\big)~(J^k_{\rho}(\Omega^{-1})^{\rho\sigma})e^{\frac12J\Omega^{-1}tJ}\Big|_{J=0}\bigg)~.\eea
To proceed further, we observe that $R^q_k$ is an operator that renames the formal variables $t_{ij}$ and does not
 touch any other index on $J$. Its effect on the last exponential is
\bea  &&R^{i_p}_k\exp\Big(\frac12\sum_{r,s}J^r_{\mu}(\Omega^{-1})^{\mu\nu}t_{rs}J^s_{\nu}\Big)\nn\\
&=&\exp\Big(\frac12\sum_{r,s\neq i_p}J^r_{\mu}(\Omega^{-1})^{\mu\nu}t_{rs}J^s_{\nu}
+\sum_{r\neq i_p}J^r_{\mu}(\Omega^{-1})^{\mu\nu}t_{rk}J^{i_p}_{\nu}\Big)\nn\\
&=&\exp\Big(\frac12\sum_{r,s\neq i_p,k}J^r_{\mu}(\Omega^{-1})^{\mu\nu}t_{rs}J^s_{\nu}
+\sum_{r\neq i_p,k}J^r_{\mu}(\Omega^{-1})^{\mu\nu}t_{rk}(J^{i_p}_{\nu}+J^{k}_{\nu})\Big)~.\nn\eea
The last relation implies the derivation of the exponential with respect to  $J^{i_p}$ is the same as
 derivation with respect to $J^k$, with this we get
\bea (\ref{stella})&=&-\sum_{k=1}^N\sum_{p=1}^l\sum_{i_1\cdots i_l=1}^N\frac{\partial}{\partial t_{i_p}}\bigg(\big(t_{i_1}f_1(\frac{\partial}{\partial J^{i_1}})\big)\cdots
\Big(t_{i_p}\partial_{\sigma}f_p(\frac{\partial}{\partial J^{k}})\Big)\cdots\nn\\
&&\hspace{3cm}\cdots \big(t_{i_l}f_l(\frac{\partial}{\partial J^{i_l}})\big)~(J^k_{\rho}(\Omega^{-1})^{\rho\sigma})e^{\frac12J\Omega^{-1}tJ}\Big|_{J=0}\bigg)~.\nn\eea
The summation over $i_p$ is now almost trivial, it can only take $N-l+1$ values because the $t_i$'s are odd. So the summation over $i_p$ gives $(N-l+1)$ times the following
\bea &&\sum_{k=1}^N\sum_{p=1}^l(-1)^p\bigg(\sum_{i_1\cdots\hat{i_p}\cdots i_l=1}^N\big(t_{i_1}f_1(\frac{\partial}{\partial J^{i_1}})\big)\cdots
\Big(\partial_{\sigma}f_p(\frac{\partial}{\partial J^{k}})\Big)\cdots\nn\\
&&\hspace{3cm}\cdots
\big(t_{i_l}f_l(\frac{\partial}{\partial J^{i_l}})\big)~(J^k_{\rho}(\Omega^{-1})^{\rho\sigma})e^{\frac12J\Omega^{-1}tJ}\Big|_{J=0}\bigg)~.\nn\eea
This is exactly the same as (\ref{star}) (upon switching $\rho,\sigma$). To conclude, we have shown
\bea \bra \partial_{CE}(f_1,\cdots f_l)\ket=-\frac{1}{2(N-l+1)}\sum_{k,q=1}^NR^q_k\frac{\partial}{\partial t_{kq}}\frac{\partial}{\partial t_{q}}\bra (f_1,\cdots f_l)\ket~.\nn\eea
Note that in the above sum, one must take the derivative $\partial_{t_{kq}}$ first then set $q=k$.
This is exactly what we are after, the integral gives the graph corresponding to ${\cal O}[f_1] ... {\cal O}[f_l]$, and $\partial_{t_{kq}}$
removes one edge from vertex $k$ to $q$ and renames the new vertex as $k$ inheriting all the other edges\footnote{If there are originally more than one edge between $j$ and $k$ then we get zero since $t_{kk}=0$. While on the graph side, contracting two vertices with more than one edge in between will give loops which also leads to zero.}. The sign factor $-(-1)^{p-1}$ is as given in subsection \ref{BOGT}. We have now arrived at a neat formula
\bea &&\bra\partial_{CE}(f_1,\cdots f_l)\ket=\partial_{Gph}\bra(f_1,\cdots f_l)\ket~.\label{iso_CE_gph-more}\eea
Thus we have completed the proof of the theorem by Kontsevich: there is a homomorphism between the Chevalley-Eilenberg complex ${\rm CE}_{\bullet}(\textrm{Ham}^0(\BB{R}^{2n}))$  and the graph complex given by the path integral. The proof of the homomorphism is the hard part of the theorem, while to prove that the mapping is in fact bijective is rather trivial. Suppose the dimension of the target space $\BB{R}^{2n}$ is big enough, one can then always find a set of polynomials, which upon applying the Feynman rules gives any given graph, this is left as an exercise for the reader.

\subsubsection{Example}

One can easily recast the example from subsection \ref{example-KT} into the calculation of the correlator  $\bra (e,f,g,h) \ket$
  according the formula (\ref{inegarakkee22}). Now we deal with graph chains as formal polynomials.

 It is instructive to show how the graph differential acts on the polynomials.
 Let us see the example of applying the differential operator (\ref{gph_diff_poly}) to some specific polynomials
\bea&&\Gamma_1(t)=\sum_{i_1,...,i_4}t_{i_1}t_{i_2}t_{i_3}t_{i_4}t_{i_1i_2}t_{i_1i_4}t_{i_1i_3}t_{i_2i_4}t_{i_4i_3}t_{i_2i_3}~,\nn\\
&&\Gamma_2(t)=\sum_{i_1,...,i_4}t_{i_1}t_{i_2}t_{i_3}t_{i_4}(t_{i_1i_2})^2(t_{i_4i_3})^2t_{i_1i_4}t_{i_2i_3}~,\nn\eea
that correspond to $\Gamma_1,~\Gamma_2$ on Figure \ref{ex_box_fig} (here these graphs are understood as equivalence classes of
 the graphs).
  A direct calculation gives (we use $i,j,k,l$ to denote $i_1,i_2,i_3,i_4$)
\bea-2\partial\Gamma_1(t)&=&\sum_{jkl}
t_{j}t_{k} t_{l}\big(t_{jl} t_{jk} t_{jl} t_{lk} t_{jk}+ t_{kj} t_{il} t_{jl} t_{lk} t_{jk}+t_{lj} t_{lk} t_{jl} t_{lk} t_{jk}\big)\nn\\
&&-\sum_{ikl} t_{i}t_{k} t_{l} \big(-t_{il} t_{ik} t_{il} t_{lk} t_{jk}+ t_{ik} t_{il} t_{ik} t_{kl} t_{lk}+ t_{il} t_{il} t_{ik} t_{lk} t_{lk}\big)\nn\\
&&+\sum_{ijl}t_{i}t_{j}t_{l} \big(-t_{ij} t_{il} t_{jl} t_{li} t_{ji} -t_{ij} t_{il} t_{ij} t_{jl} t_{lj}- t_{ij} t_{il} t_{il} t_{jl} t_{jl}\big)\nn\\
&&-\sum_{ijk}t_{i}t_{j}t_{k} \big(-t_{ij} t_{ik} t_{ji} t_{ik} t_{jk}- t_{ij} t_{ij} t_{ik} t_{jk} t_{jk} +t_{ij} t_{ik} t_{ik} t_{jk} t_{jk}\big)\nn\\
&=&-12\sum_{ijk}t_{i}t_{j}t_{k}t^2_{ij} t^2_{ik} t_{jk}~,\nn\\
-2\partial\Gamma_2(t)&=&
\sum_{jkl}t_{j}t_{k} t_{l} \big( t_{lj} ^2 t_{lk} ^2 t_{jk}\big)-\sum_{ikl}t_{i}t_{k} t_{l} \big(-t_{ik} ^2 t_{lk} ^2 t_{il}\big)\nn\\
&&+ \sum_{ijl}t_{i}t_{j}t_{l}\big(-t_{ij} ^2 t_{lj} ^2 t_{il}\big)-\sum_{ijk} t_{i}t_{j}t_{k} \big(-t_{ij} ^2 t_{ik} ^2 t_{jk}\big)\nn\\
&=&4\sum_{ijk}t_{i}t_{j}t_{k} t_{ij} ^2 t_{ik} ^2 t_{jk}~.\nn\eea
To extract the coefficient of $\Gamma_8$ from these polynomials, one applies the operator
\bea \Gamma_8^*=\frac18\frac{\partial}{\partial t_3}\frac{\partial}{\partial t_2}\frac{\partial}{\partial t_1}\frac{\partial^2}{\partial t_{12}^2}\frac{\partial^2}{\partial t_{13}^2}\frac{\partial}{\partial t_{23}}~,\nn\eea
  which can be understood as graph cochain (see Appendix \ref{graph-cocycle-proof}).
For example
\bea &&\bra\Gamma^*_8, \partial\Gamma_1\ket=\frac68\frac{\partial}{\partial t_3}\frac{\partial}{\partial t_2}\frac{\partial}{\partial t_1}\frac{\partial^2}{\partial t_{12}^2}\frac{\partial^2}{\partial t_{13}^2}\frac{\partial}{\partial t_{23}}\sum_{ijk=1}^Nt_{i}t_{j}t_{k}t^2_{ij} t^2_{ik} t_{jk}\nn\\
&=&\frac34\frac{\partial^2}{\partial t_{12}^2}\frac{\partial^2}{\partial t_{13}^2}\frac{\partial}{\partial t_{23}}\Big(t^2_{12} t^2_{13} t_{23}-t^2_{13} t^2_{12} t_{32}-t^2_{21} t^2_{23} t_{13}+t^2_{23} t^2_{21} t_{31}+t^2_{31} t^2_{32} t_{12}-t^2_{32} t^2_{31} t_{21}\Big)\nn\\
&=&6~.\nn\eea
We end up with the same results as before. This is just a simple illustration how to work with the graphs as formal polynomials. We do want to point out that for concrete graphs it is much more efficient to work with pictures than to use the polynomial representation. However, the latter approach allows us to prove theorems about graph complex by manipulating operators acting on the polynomials;  thus
  instead of drawing lots of pictures one 'lets algebra do the talking'.

\subsection{Generalizations}

Let us briefly sketch the possible generalizations of the above construction about the relation between certain
 Chevalley-Eilenberg complex and certain graph complex. As the reader may imagine there are many ways to generalize and
  extend the presented construction. Here we just indicate some directions and leave to the reader to figure out the details.

 The most straightforward generalization is related to case of $\BB{R}^{2n|m}$ equipped with even symplectic structure
 \bea
  \Omega_{\mu\nu}~ dx^\mu \wedge dx^\nu + \eta_{ab} ~d\psi^a \wedge d\psi^b~,\nn
 \eea
  where $\Omega$ is non-degenerate $2n\times 2n$ antisymmetric matrix and $\eta$ is symmetric $m \times m$ matrix.
 Here $x^\mu$ ($\mu=1,2,..., 2n$) stands for even coordinate of  $\BB{R}^{2n|m}$ and $\psi^{a}$ ($a=1,2,...,m$) for the odd
  coordinate of $\BB{R}^{2n|m}$. In analogy with our previous discussion we can define the Chevalley-Eilenberg complex
 ${\rm CE}_{\bullet}(\textrm{Ham}^0(\BB{R}^{2n|m}))$ for the formal Hamiltonian vector fields on $\BB{R}^{2n|m}$.
  In analogy with  Theorem \ref{theorem-cocycle}  we can formulate the following theorem on the relation between  the Chevalley-Eilenberg complex
 ${\rm CE}_{\bullet}(\textrm{Ham}^0(\BB{R}^{2n|m}))$  and the graph complex $\Gamma_\bullet$.

    \begin{Thm}[super-Kontsevich]\label{super-theorem-cocycle}
    Let us take a collection of polynomial functions $f_1(x, \psi)$,  ..., $f_l(x, \psi)$ on $\BB{R}^{2n|m}$ and define the following integral
\bea \bra (f_1, ... , f_l) \ket=\frac{1}{Z[0]} \int\limits_{\bigoplus\limits_{i=1}^N \mathbb{R}^{2n|m}} {\cal D}x {\cal D} \psi
~\sum_{i_1=1}^N t_{i_1}f_1(x_{i_1}, \psi_{i_1}) ... \sum_{i_l=1}^N t_{i_l}f_l(x_{i_l}, \psi_{i_l})~e^{-S}~\in \BB{R}[t_i,t_{ij}] \label{super-inegarakkee22}\eea
where ${\cal D}x= d^{2n}x_1\cdots d^{2n}x_N$ and ${\cal D} \psi = d^m\psi_1 \cdots d^m \psi_m$ and
\bea
S = \frac12\sum\limits_{i,j=1}^N t^{ij} \left (  x_i^\mu  x_j^\nu ~\Omega_{\mu\nu} + \psi^a_i \psi^b_j \eta_{ab} \right ) ~.\label{super-BV_action_discrete}\eea
  The integral (\ref{super-inegarakkee22})
 can be understood as map from chain $(f_1, ... , f_l)$ in ${\rm CE}_{\bullet}(\textrm{Ham}^0(\BB{R}^{2n|m}))$ to graph chain
  in $\Gamma_\bullet$.  Moreover the map (\ref{inegarakkee22}) satisfies the property
\bea
\bra\partial (f_1,\cdots f_l)\ket=\partial_{Gph}\bra(f_1,\cdots f_l)\ket~,\label{iso_CE_gph}
\eea
where $\partial$ is the appropriate Chevalley-Eilenberg differential and $\partial_{Gph}$ is the graph differential defined in (\ref{gph_diff_poly}).
\end{Thm}

To understand (\ref{super-inegarakkee22}) as concrete prescription we have to develop and study the perturbative expansion in the odd-coordinates.  The proof of this theorem is straightforward generalization of the proof for Theorem \ref{theorem-cocycle}. We have to embed
 the integral (\ref{super-inegarakkee22})  on $\bigoplus\limits_{i=1}^N \mathbb{R}^{2n|m}$
  to BV theory on $\bigoplus\limits_{i=1}^N \mathbb{R}^{2n|m} \oplus  \bigoplus\limits_{i=1}^N \mathbb{R}^{2n|m}[-1]$ and use the Ward identities. Here the only complication compared with the previous discussion is
  related to some additional signs due to presence of odd coordinates $\psi$'s.  We leave all these details to the reader to figure out.

 By looking at Theorems \ref{theorem-cocycle} and \ref{super-theorem-cocycle} it is very easy to construct the graph
   cycles by choosing simple Chevalley-Eilenberg
  cycles.  Let us give three different examples.

  \begin{Exa}\label{pure-odd-cocycle}
   Consider the case $\BB{R}^{0|m}$ and let us choose the cubic function
   $$f(\psi) = \frac{1}{3!} f_{abc} \psi^a \psi^b \psi^c$$
 such that $\{f, f\} =0$.  This is equivalent to the statement that $f^{a}_{~bc}= \eta^{ad} f_{dbc}$ are structure constants
  for some Lie algebra $\FR{g}$ and $\eta_{ab}$ is $ad$-invariant metric on this Lie algebra (here $\eta^{ad} \eta_{db} = \delta^a_b$).
  Using the prescription (\ref{super-inegarakkee22}) we can construct the graph cycle
   since $(f,f,...,f)$ is trivially Chevalley-Eilenberg cycle due to property $\{f, f\}=0$.  The corresponding
    cycle is constructed as
    \bea
     \sum\limits_{\Gamma} c_{\Gamma} (\FR{g})~ \Gamma~,
    \eea
   where $\Gamma$'s are collection of trivalent graphs and $c_\Gamma$ is the number constructed by contraction
    of $f_{abc}$ as vertices and $\eta^{ab}$ as propagator according to the graph $\Gamma$.
 The numbers $c_{\Gamma}$ satisfies so-called  IHX-relations (see \cite{BN1},\cite{sawon1999}) which is equivalent to the statement that above
  expression is a graph cycle.
\end{Exa}

\begin{Exa}\label{pure-even-cocycle}
 Consider the standard symplectic vector space $\BB{R}^{2n|0}$ and let us choose the collection of function $f_i$ ($i=1,2,..., l$)
  such that $\{f_i, f_j\}=0$. For example, the case $l=n$ would correspond to the completely integrable system on $\BB{R}^{2n}$.
   Plug these collection of functions into the prescription of Theorem (\ref{theorem-cocycle}) and we will get the graph cycle since
    $(f_1, f_2, ... , f_l)$ is trivially Chevalley-Eilenberg cycle due to property $\{f_i, f_j\}=0$.  The corresponding cycle
     has the form
     \bea
      \sum\limits_{\Gamma} c_\Gamma (f_1, ..., f_l)~ \Gamma~,
     \eea
   where  $c_\Gamma (f_1, ..., f_l)$ are numbers constructed from the contraction of vertices (Taylor coefficients of $f$'s)
      using the inverse of symplectic structure $\Omega^{-1}$ according to the graph $\Gamma$.  The concrete prescription with
       all numerical coefficients can be read off the perturbative expansion of (\ref{inegarakkee22}).
\end{Exa}

\begin{Exa}\label{Linfty-cocycle}
 Now let us consider the general case $\BB{R}^{2n|m}$ and combine these two examples above. First of all
  we can find an odd function $\Theta(x, \psi)$ such that $\{\Theta, \Theta \}=0$. Indeed such $\Theta$  is related
   to $L_\infty$-algebra with invariant metric (here the symplectic structure on $\BB{R}^{2n|m}$ is such a metric).
    The corresponding graph cycle is constructed by contracting the Taylor coefficients of $\Theta$ as vertices
    and the inverse of symplectic structure as propagators. The precise prescription can be read off from the perturbative
     expansion of (\ref{super-inegarakkee22}) which would allow as to construct graph cycle starting from a cyclic $L_\infty$-algebra.

     Also in analogy with example \ref{pure-even-cocycle} we can choose the collection of $\Theta_i(x,\psi)$, $i=1,2,...,l$
      such that $\{ \Theta_i, \Theta_j \}=0$. Plugging these $\Theta$'s into (\ref{super-inegarakkee22}) we will get the graph cycle.
\end{Exa}

We can continue to play this game and construct more exotic BV theories and get more exotic graph complexes. For example,
 for the graph complex related to chord diagrams and the related Chevalley-Eilenberg complex see \cite{2010arXiv1006.1240Q}.

\section{Outline for quantum field theory}
 \label{QFT-summary}

So far our discussion involved the finite dimensional integrals and the related algebraic structures. In quantum field theory
 we have to study the infinite dimensional integrals and construct the corresponding BV formalism with the Ward identities.
  Referring to the terminology of section \ref{oddFT} now the manifold $M$ is infinite dimensional functional space of fields,
   the derivatives should be replaced by functional derivatives and formally the odd cotangent bundle $T^*[-1]M$ should
    be  BV manifold.
    However we suffer from the standard problems with infinite dimensional setting, namely some of the formulas are not well-defined
     and may require the additional regularization.  For example, there is no simple canonical way to define the odd Lapalacian (\ref{def-Laplace})
      in functional space since the double functional derivative is not well-defined as it stands.  At the same time the odd
       Poisson bracket (\ref{def-PBoncot}) can be defined in suitable way in infinite dimensional setting.  Nevertheless one tries to
        proceed formally and apply the Ward identities (\ref{Ward1}) and (\ref{Ward2}) formally. Despite all these problems
         the BV formalism provides good heuristic understanding of infinite dimensional path integral.
 The  impressive example of using BV formalism in infinite dimensional setting is given by the
        heuristic derivation of deformation quantization and the proof of formality in the context of Poisson sigma model, \cite{Cattaneo:1999fm}.

   However we will not discuss the infinite dimensional BV formalism here.
  We restrict ourselves to a few general remark regarding the perturbative aspects of QFT. In particular we would like to
  concentrate on the Chern-Simons type of theories and their different relatives.  We will introduce the theory in purely combinatorial fashion
 without any reference to infinite dimensional formulation. Since any perturbative expansion of QFT is done through the graphs
  (Feynman diagrams), then we would like to define certain nice objects on graphs such as element of graph homology and
   element of graph cohomology. The infinite dimensional path integral (partition function) will be defined as pairing of these
    two elements. We believe that this picture is rather generic for the perturbative expansions of QFT. However the details has been
     worked out only for some particular examples.  Let us illustrate briefly the idea on the example of Chern-Simons theory and its different
      modification.

   \subsection{Formal Chern-Simons theory and graph cocycles}

    In previous section we have described how one can construct the graph cycles, see examples \ref{pure-odd-cocycle},
    \ref{pure-even-cocycle} and \ref{Linfty-cocycle}.
  At the same time
    the QFT theory gives us a very natural and elegant way to construct the cocycles. Let us illustrate this by the example inspired by
    the Chern-Simons theory on $S^3$ (or any rational homology sphere where $H^1(M, \BB{R})=0$).  Let us consider
     the supermanifold $T[1]S^3$ with coordinates $z = (\theta^a, \sigma^a)~,a=1,2,3$ and define homological vector field
 $D_z = \theta^a \frac{\partial}{\partial \sigma^a}$. We can introduce the propagator $G(z_1, z_2)$ which is smooth function defined on
 $(z_1, z_2) \in T[1]S^3 \times T[1]S^3$ minus the diagonal ($z_1 \neq z_2$). The propagator satisfies the following equation
\bea
\left ( D_{z_1} + D_{z_2} \right ) G(z_1, z_2) = -\delta^6 (z_1 - z_2)~,\label{defin-propagator-CS}
\eea
 where  $\delta^6 (z)$ is delta-function on $T[1]S^3$ with canonical integration.  In certain sense the propagator is inverse of
  de Rham operator $D_z$ and  obviously it is not uniquely defined. $D_z$ can be inverted on co-exact forms and thus
   would require the Hodge decomposition with a concrete metric. The ambiguity in (\ref{defin-propagator-CS}) is given
    by the following formula
\bea
 G(z_1, z_2)~ \rightarrow~ G(z_1, z_2) + \left ( D_{z_1} + D_{z_2} \right ) L(z_1, z_2)~,\label{change-of-propogator}
\eea
 where $L(z_1, z_2)$ is some function on $T[1]S^3 \times T[1]S^3$.  Assume that we choose some $G(z_1, z_2)$ satisfying
  the above properties, then we can construct the following differential operator acting on the graph chain $\Gamma(t)$
\bea b_{\Gamma}=\exp\Big(\sum^N_{i=1}\int  d^6z_i \frac{\partial}{\partial t_i} \Big)~\exp \Big( \sum^N_{k,l=1} G(z_k,z_l) \frac{\partial}{\partial t_{kl}}\Big)~\Gamma(t)\Big|_{t_{\cdot\cdot}=0}~, \label{defin-cohain-CS}
\eea
 where $\Gamma(t)$ is understood as polynomials in $t_i$ and $t_{ij}$ (see subsection \ref{ADoGC}).  Thus $b_\Gamma$ gives a number
   for a given graph chain $\Gamma(t)$.  In  $b_{\Gamma}$ we replace every edge  $t_{ij}$ with an actual propagator running
   from $z_i$ to $z_j$ and every vertex $t_i$ with the integration $\int d^6 z_i = \int d^3\theta_i d^3\sigma_i$  over $T[1]S^3$.
The reader  may check that in $b_\Gamma$ both symmetry factors $\#V,\#P$ of a graph are taken care of automatically.
 Thus the prescription (\ref{defin-cohain-CS}) can be understood as way of producing the graph cochain and can be symbolically
  written as
 \bea
  b_{\Gamma} = \Gamma^* \left ( \frac{\partial}{\partial t} \right ) \Gamma(t) \Big|_{t_{\cdot\cdot}=0}~.
 \eea

One can prove the following statements about this cochain
\begin{itemize}
\item  $b_\Gamma$ is well-defined;
\item $b_{\partial_{Gph} \Gamma} = \delta_{Grp} b_{\Gamma}=0$ and $b_\Gamma$ is graph cocycle;
\item  under change (\ref{change-of-propogator}) of the propogator $G$, $b_\Gamma$ changes by coboundary $b_\Gamma \rightarrow b_\Gamma + \delta_{Grp} (...)$;
\end{itemize}

Thus $b_\Gamma$ as an element of graph cohomology is well-defined and is independent of the choice of the correlator.
 The proof of the first statement about $b_\Gamma$ being well-defined can be found in \cite{Axelrod:1991vq}. Two other statements are  known,
  however we cannot find a complete proof of them in the literature. Thus we provide a proof of these two statement
   in Appendix \ref{graph-cocycle-proof} for any differential graded Frobenius algebra with trivial cohomology (maybe except the top degrees).
    The space $C^\infty (T[1]S^3)$ is an example of such infinite dimensional algebra. As long as all expressions are well-defined
 the proofs for finite  and infinite dimensional cases are identical.

 Next we can take the statement of Theorems \ref{theorem-cocycle} and \ref{super-theorem-cocycle} about isomorphism of
  the Chevalley-Eilenberg  complex and graph complex and by applying the differential operator (\ref{defin-cohain-CS}) to
   (\ref{inegarakkee22}) (or (\ref{super-inegarakkee22}))
   we construct the Chevalley-Eilenberg cochain
   \bea
    c^l (f_1, f_2, ... , f_l) \equiv  \Gamma^*\left (\frac{\partial}{\partial t} \right )
     \bra (f_1, f_2, ... , f_l )\ket  \Big|_{t_{\cdot\cdot}=0}~.\label{dsjdo3o3odjs}
   \eea
    As a simple consequence of these theorems and the fact that $b_\Gamma$ is cocycle, we conclude that $c^l$ is
     the Chevalley-Eilenberg cocycle. Moreover under the change of the propagator (\ref{change-of-propogator}) this cocycle will change by
    a coboundary.  Thus $c^l$ as an element of  the Chevalley-Eilenberg cohomology depends only on the concrete Frobenius algebra
       (here related to $T[1]S^3$).
       This is an example of the understanding of the infinite dimensional integral as a coycle with certain specific properties.
        Indeed the expression (\ref{dsjdo3o3odjs}) can be represented as an infinite dimensional integral and the properties stated in
         Theorem \ref{theorem-cocycle} can be derived by some formal manipulations with this integral.

Moreover we can define the partition function as follows.
 Once there is a graph cocycle we can construct the Chern-Simons partition function for the Lie algebra $\FR{g}$ by pairing the cocycle
  (\ref{defin-cohain-CS}) with graph cycle constructed in the example \ref{pure-odd-cocycle}
  \bea
   Z[\FR{g}] = \sum\limits_{\Gamma} b_\Gamma c_{\Gamma} (\FR{g}) = \sum\limits_\Gamma \Gamma^* \left ( \frac{\partial}{\partial t} \right ) c_{\Gamma} (\FR{g}) \Gamma(t) \Big|_{t_{\cdot\cdot}=0}~.
  \eea
   Since $\sum\limits_{\Gamma} c_{\Gamma} (\FR{g}) \Gamma$ is graph cycle,  any changes of $b_\Gamma$ by a coboundary
    vanish upon the pairing. Thus the expression   $Z[\FR{g}]$  is independent of the concrete choice of the propagator and it depends only
     on Lie algebra $\FR{g}$ and the Frobenius algebra ($T[1]S^3$ in the case of standard Chern-Simons theory). Thus in this
      sense  $Z[\FR{g}]$ can be though as an invariant of $S^3$ for fixed $\FR{g}$.

  Analogously we can pair the cocycle (\ref{defin-cohain-CS}) with any other graph cycles, see example \ref{pure-even-cocycle} and \ref{Linfty-cocycle}.  These will give rise to some different partition functions. For instance, the example \ref{Linfty-cocycle}
   produces the partition function depending on cyclic $L_\infty$-algebra and Frobenius algebra (for example, $T[1]S^3$ as done in \cite{QiuZabzine:2009rf}).
    This sort of models were discussed in \cite{QiuZabzine:2009rf} and they are natural generalizations of Chern-Simons theory.

 Moreover instead of $T[1]S^3$ one may use any acyclic Frobenius algebra for the construction of cocycles and thus
  deal with the formal and discrete versions of Chern-Simons theory. For the examples of
 discrete version of Chern-Simons the reader may consult  \cite{Cattaneo:2008ph}.

 Depending of concrete QFT and the set of observables  we may be forced to study more complicated graph complexes which can
  be colored, decorated, have external legs or other additional structures. However we would expect that the qualitative picture
   of perturbation theory as a pairing of appropriate graph cycle with graph cocycle will still hold good. These issues would require the additional
    study.

\bigskip\bigskip

\noindent{\bf\Large Acknowledgement}:
\bigskip

  M.Z. is grateful to the organizers of the Winter School
 for the invitation and for the warm hospitality during his stay in Czech Republic.
  We thank Rikard von Unge
    for the reading and commenting on these lecture notes.
The research of M.Z. is supported by  VR-grant 621-2008-4273.

\bigskip\bigskip

\appendix
\section{Explicit formulas for odd Fourier transform}\label{app-A}

This appendix should be regarded as companion for subsection \ref{ss-odd-FT} and we follow the same
 notations as in  subsection \ref{ss-odd-FT}. We give here some explicit formulas and derive some curious relations.

  The odd Fourier transform (\ref{F-T}) maps the differential forms to multivectors according the following
   explicit formula
\bea
&& \frac{1}{p!} f_{\mu_1 ... \mu_p} (x) dx^{\mu_1} \wedge ... \wedge dx^{\mu_p}
  ~\overset{F}{\longrightarrow}~ \nn\\
&&\hspace{3.5cm}\frac{(-1)^{(n-p)(n-p+1)/2}}{p! (n-p)!} f_{\mu_1... \mu_p} \Omega^{\mu_1... \mu_p \mu_{p+1} ... \mu_n} \partial_{\mu_{p+1}} \wedge
 ... \wedge \partial_{\mu_n}\label{FT-explicit}~,
\eea
  where $\Omega^{\mu_1... \mu_n}$ is defined as components of a nowhere vanishing
   top multivector field dual to the volume form (\ref{defin-volume})
  $$ {\rm vol}^{-1} = \rho^{-1}(x)~ \partial_1 \wedge ... \wedge \partial_n= \frac{1}{n!}~ \Omega^{\mu_1 ... \mu_n}(x)~\partial_{\mu_1} \wedge...
 \wedge \partial_{\mu_n}~.$$
 The operator $\Delta$ corresponds to the divergence operator acting on the multivector fields as
  follows
  $$  ({\rm div} ~\tilde{f})^{\mu_1 ... \mu_p} \partial_{\mu_1} \wedge ... \wedge \partial_{\mu_{p-1}} =
   \frac{1}{(p-1)!} ~ \rho^{-1} \frac{\partial}{\partial x^\nu} \left ( \rho \tilde{f}^{\nu \mu_1 ... \mu_{p-1}} \right )
    \partial_{\mu_1} \wedge ... \wedge \partial_{\mu_{p-1}}~,$$
     where we use the convention for the identification of the multivector with function on the odd cotangent bundle
      from the example \ref{oddcotangentbundle}.  The odd Poisson bracket (\ref{def-PBoncot}) corresponds to
       the Schouten bracket
\bea
&&\left (\{ \tilde{f}, \tilde{g} \}\right )^{\mu_1 ... \mu_{p+q-1}} \partial_{\mu_1} \wedge ... \wedge \partial_{\mu_{p+q-1}} \nn \\
&&=\frac{1}{p! q!} \left (q~ \partial_\mu f^{\nu_1... \nu_p} g^{\mu \nu_p ... \nu_{p+q-1}} + (-1)^p p ~f^{\mu \nu_1 ... \nu_{p-1}}
 \partial_\mu g^{\nu_p ... \nu_{p+q-1}} \right)  \partial_{\mu_1} \wedge ... \wedge \partial_{\mu_{p+q-1}}~, \nn
 \eea
  which is the generalization of the usual Lie bracket to the multivector fields.

  Let us mention a few curious facts about the transportation of the BV algebra structure $(C^\infty (T^*[-1]M), \cdot, \{~,~\}, \Delta)$ to $C^\infty (T[1]M)= \Omega^\bullet (M)$ using the inverse Fourier transform (\ref{inverse-FT}).
   The graded commutative product on  $C^\infty (T^*[-1]M)$ get mapped to the following product on
    $C^\infty (T[1]M)$
    \bea
     f * g = F^{-1} \left ( F[f] F[g] \right ) = (-1)^{(|f|+n)n} \int d^n \xi~ \rho^{-1}~ f(x, \xi) g(x, \theta - \xi)~,\label{newstartdifff}
    \eea
     where $\xi$ and $\theta$ are odd coordinates on $T[1]M$.  This star  product is associative product of degree
      $-n$ and thus $|f*g| = |f| +|g| -n$.  The commutativity rule is $f*g = (-1)^{(n-|f|)(n-|g|)} g*f$ and thus in
       general it is not a graded commutative product. If we work with $\mathbb{Z}_2$-grading and $n=\dim M$ is even,
        the product (\ref{newstartdifff}) is supercommutative. The odd Poisson bracket (\ref{def-PBoncot}) gives rise to the
         bracket on $C^\infty (T[1]M)$
         \bea
     [f, g] = F^{-1} \left ( \{ F[f], F[g] \} \right ) =   (-1)^{|f|} D ( f* g) - (-1)^{|f|} (Df) * g - (-1)^n f * (Dg)~,
         \eea
  where to derive the last relation  we used (\ref{D-Delta-relation}), (\ref{Delta-defPB}) and (\ref{newstartdifff}).
   The bracket $[~,~]$ is  of degree $(1-n)$.
   One can easily derive the properties of the bracket $[~,~]$ on $C^\infty (T[1]M)$ by the Fourier transform
    of the properties of the Gerstenhaber algebra (see the definition \ref{Gerstenhaber-algebra}).  Let us just
     point out that in the case of $\mathbb{Z}_2$-grading and $n=\dim M$ being even the bracket $[~,~]$ is
      Gerstenhaber bracket with respect to the supercommutative multiplication $*$.

\section{BV-algebra on differential forms}\label{app-B}

 If the manifold $M$ is equipped with a Riemannian metric $g_{\mu\nu}$ then we can define the odd Fourier
  transform which maps $C^\infty (T[1]M)$ to $C^\infty (T[-1]M)$. Using the metric we can define the odd Fourier transform
for $f(x, \theta) \in C^\infty (T[1]M)$ as follows
\bea
 F[f](x, \xi) = \int d^n\theta~g^{-1/2} e^{\xi^\mu g_{\mu\nu} \theta^\mu}  f(x, \theta)~,\label{F-T-diff.forms}
\eea
 where $g = \det (g_{\mu\nu})$ and $\xi^\mu$ is odd coordinate of degree $-1$ on $T[-1]M$. If we work in
  superlanguage then the odd Fourier transform (\ref{F-T-diff.forms}) maps $C^\infty (\Pi TM)$ to itself.
 In the language of differential forms the odd Fourier transform corresponds to Hodge star operation $\star$
\bea
 F[f_{(p)}] = (-1)^{(n-p)(n-p+1)/2}*f_{(p)}~,\nn\eea
  where $f_{(p)}$ is $p$-form.  Under the Fourier  transform the homological vector field $D$ is mapped as
 follows
$$D^\dagger F[f] = (-1)^n F[Df]~,$$
 where $D^\dagger$ is defined as
 \bea
  D^\dagger = \frac{1}{\sqrt{g}} \frac{\partial}{\partial x^\mu} g^{\mu\nu}(x) \frac{\partial}{\partial \xi^\nu} \sqrt{g}
   - \Gamma^\gamma_{~\sigma\mu} g^{\sigma \nu} \xi^\mu \frac{\partial}{\partial \xi^\gamma} \frac{\partial}{\partial
    \xi^\nu}~,
 \eea
  where $ \Gamma^\gamma_{~\sigma\mu}$ is the Levi-Civita connection for the metric $g_{\mu\nu}$.
  $D^\dagger$  is operation of degree $1$ on $C^\infty(T[-1]M)$ and by construction $(D^\dagger)^2=0$. On the space of
   differential forms $\Omega^\bullet (M) = C^\infty(T[-1]M)$ the operator $D^\dagger$ is proportional
    to the adjoint of the de Rham differential, $d^\dagger$.  The operator $D^\dagger$ is a second order differential operator and it satisfies the relation
   (\ref{Delta-quadratic}) with the usual graded commutative multiplication on $C^\infty(T[-1]M)$.
    Therefore the space of differential forms $\Omega^\bullet (M)$ is equipped with the BV-structure where $D^\dagger$
     corresponds to $\Delta$ and the odd bracket is
     $$ \{f, g\} =  (-1)^{|f|} D^\dagger (fg)+(-1)^{|f|+1}(D^\dagger f)g - f(D^\dagger g)  $$
  where $f,g \in \Omega^\bullet (M)$. The integration theory is canonically defined for the differential form
   and thus it is not hard to extend the discussion from subsections \ref{s-integration} and \ref{AVOTI} to the present
    case of differential forms.

  The story presented here can be reiterated in many other cases, maybe with some minor modifications.
   For example, instead of metric we can use   the symplectic structure in the definition of the odd Fourier transform
   (\ref{F-T-diff.forms}).  We can also treat   in similar fashion the Lie  algebroid which corresponds to a vector bundle with odd fiber coordinate  $A[1]$ and with the homological field $Q$ of degree $1$.
    The integration on $A[1]$ will require some additional structure. The odd Fourier transform will map
     $C^\infty (A[1])$ to $C^\infty (A^*[-1])$ etc.
     Thus we can get numerous examples of BV algebras and related structures.

\section{Graph Cochain Complex}
\label{graph-cocycle-proof}

The graph cochain complex is defined as a formal linear combination of dual graphs
\bea \Gamma^*=\sum_{\G}b_{\G}\G^*~,\nn\eea
and the evaluation of the graph cochain on a graph chain is often written as a pairing
\bea \bra \sum_{\G}b_{\G}\G^*,\sum_{\G}c_{\G}\G\ket\in \BB{R}\textrm{ or }\BB{C}~.\nn\eea
The differential on co-graphs are defined as
\bea \bra \delta_{Gph}\sum_{\G}b_{\G}\G^*,\sum_{\G}c_{\G}\G\ket=\bra\sum_{\G}b_{\G}\G^*,\partial_{Gph}\sum_{\G}c_{\G}\G\ket~.\nn\eea

Like the graph chain complex, the graph chain complex can be presented as polynomials of differential operators in $t_{ij},t_i$
\bea \Gamma^*(\frac{\partial}{\partial t_{ij}},\frac{\partial}{\partial t_{i}})~,\nn\eea
which acts on the graph polynomials in the obvious way. Instead of writing $\partial_t$, we will just name some new formal 'dual momentum' $s^{ij},s^i$ of $t_{ij},t_{i}$, and write the polynomials of differentials as polynomials in the $s$'s.

In the graph chain complex case, we have a very neat homomorphism between certain  Chevalley-Eilenberg complex and graph complex and this homomorphism leads us to construct graph cycles from  Chevalley-Eilenberg cycles. Inspired by Chern-Simons theory
  Kontsevich gave a prescription of constructing graph cocycles from Frobenius algebra with some (rather strong) extra conditions. The idea is to construct out of the Frobenius algebra 'propagators' and 'vertex functions'; then one replaces $t_{ij},t_i$ with propagators and vertex functions.


The data needed is an acyclic differential graded Frobenius algebra $(\FR{a}, \cdot,d, \bra \cdot, \cdot \ket)$:
\begin{itemize}
\item  $\FR{a}$ is graded commutative algebra;
\item  $d$ is differential on $\FR{a}$ ($d^2=0$ and $d$ is a derivation of degree $1$);
\item  there is non-degenerate pairing $\bra\cdot,\cdot\ket\to \BB{R}$  which is compatible with the graded commutative multiplication
\bea \bra ab,c\ket=\bra a,bc\ket=\bra 1,abc\ket~;\label{pairing}\eea
\item the differential $d$ is compatible with pairing $\bra\cdot,\cdot\ket$
\bea
 \bra da, b \ket + (-1)^{|a|} \bra a, db \ket=0 ~.
\eea
\end{itemize}
Acyclic means that the cohomology of $d$ is empty, maybe except the lowest and highest degrees.

Let us give both finite and infinite dimensional examples of acyclic differential graded Frobenius algebra.
\begin{Exa}\label{example-FA-SU2}
The finite dimensional example is given by  $\FR{a}=C^{\infty}(\FR{su}(2)[1])$ with $d=f_{ab}^ce^ae^b\partial_{e^c}$,
 where $f_{ab}^c$ are the structure constants of the Lie algebra $\FR{su}(2)$.  The only non-zero scalar product is
\bea \bra e^ae^b,e^c\ket=f^{abc},\nn\eea
note that the product $e^ae^b$ is regarded as one element. This pairing is clearly graded symmetric $f^{abc}=f^{cab}$ and
 is non-degenerate for $\FR{su}(2)$ case. The cohomology of $d$ is trivial except $H^0 (d) = H^3(d) = \BB{R}$.
\end{Exa}
\begin{Exa}
Consider a rational homology $n$-sphere $\Sigma$ (i.e., $n$-manifold with the property $H^k (\Sigma, \BB{R})=0$ for $0< k < n$).
 The space of differential forms $\Omega^\bullet (\Sigma) = C^\infty (T[1]\Sigma)$ is an infinite dimensional acyclic differential Frobenius
  algebra with the multiplication given by wedge and differential by the exterior derivative.
The pairing is given by  the integral
\bea \bra \alpha ,\beta \ket=\int_{\Sigma}~ \alpha \wedge \beta~,\nn\eea
where $\alpha$ and $\beta$ are differential forms.
\end{Exa}

Let us choose the basis $e^I$ for the \emph{underlying module} of $\FR{a}$, let $f^I$ be the formal parameters corresponding to $e^I$, but with no algebraic relations besides the graded commutativity.

Because of the property (\ref{pairing}), we write the pairing suggestively as an integral
\bea \smalint e^Ie^J\cdots e^K=\bra e^I, e^J\cdots e^K\ket\in \BB{R}~.\nn\eea
We assume that $\int$ has degree $-p$, which is required to be odd (the reason will become clear later)\footnote{In principle, it is enough
 to consider $\int$ to be odd with non-homogeneous terms of different odd degrees. Although we do not know any concrete examples of this situation.}.
Define also the matrix
\bea m^{IJ}=\smalint e^Ie^J,~~\deg e^I+\deg e^J=p~,\nn\eea
and its inverse
\bea  m^{IK}m_{KJ}=\delta^I_J,~~e_I=e^Jm_{JI},~~\int e^Ie_J=\delta^I_J~.\nn\eea
We denote the differential in the matrix notation as
\bea de^I=D^I_{~J}e^J~.\nn\eea
The Stokes theorem
\bea \smalint d(\cdots)=0\nn\eea
plus the fact that $d$ is a derivation imply
\bea D^I_{~K}m^{KJ}+(-1)^Im^{IK}D^J_{~K}=0,\nn\\
m_{IL}D^L_{~J}+(-1)^{p-I}D^L_{~I}m_{LJ}=0~,\label{sym_pairing}\eea
note we use the short hand notation $(-1)^I$ where $I$ should be understood as the degree of $e^I$.

With the assumption of acyclicity, $d$ may be inverted in certain sense. In fact we shall assume that the inverse of $d$ is obtained by means of Hodge decomposition. In complete analogy with the standard Hodge theory we can pick up the metric on $\FR{a}$ and
 construct the Hodge theory for $d$.
 Thus we can Hodge decompose any element into
\bea \psi=\psi^h+d\ga+d^{\dagger}\gb~,\nn\eea
acyclicity implies the harmonic element is zero, $\psi^h=0$.
The basis $e^I$ can be chosen as the eigen-modes of the Laplacian $\square=\{d,d^{\dagger}\}$.
The inverse can be explicitly written as
\bea d^{-1}\psi=\int \psi\,K~,~~~K=e_I\otimes \frac1{\square}d^{\dagger}e^I~.\nn\eea
In the general case the propagator $K$ is written as
\bea &&K_{IJ}=m_{IK}(D^{-1})^{K}_{~J}=(-1)^{p-I}(D^{-1})^{K}_{~I}m_{KJ}~,\nn\\
&&K_{QP}=(-1)^{QP+1}K_{PQ}~.\label{sym_K}\eea
Here both properties of (\ref{sym_K}) follow from the symmetry properties of (\ref{sym_pairing}).
It shall be proved later that the details of how $d^{-1}$ is obtained does not affect the cohomology class of co-graphs.


Now we proceed to the construction of graph cocycles. First we define a formal integration operator that acts on polynomials of $f^I$
\bea \int \phi(f)=\smalint \exp\Big\{e^I\frac{\partial}{\partial f^I}\Big\}~\phi(f)\bigg|_{f=0}~.\label{integration}\eea
We can derive the following concatenation property of the integration operator. We take $N$ copies of $e_i^I,f_i^I$, and let $\int_i$
 be as in (\ref{integration}), but for the $i^{th}$ copy. Calculate the commutator
\bea
&&\Big[s^i\int\limits_is^j\int\limits_j,~ \big(m_{IJ}f^I_if^J_j(-1)^{pJ}\big)\Big]\nn\\
&=&s^j\int\limits_j~s^i\int\limits_i\exp\Big\{e_i^I\frac{\partial}{\partial f_i^I}\Big\}e_{iJ}e^J_j\exp\Big\{e_j^I\frac{\partial}{\partial f_j^I}\Big\}(-1)^{pJ}~.\label{temp1}\eea
Now notice the relation
\bea&&\Big(\int\limits_i e_i^JP_Je_{iI}\Big)e^I_j(-1)^{pI}=\Big(\int\limits_i e_i^Je_{iI}P_J(-1)^{P_J(p-I)}\Big)e^I_j(-1)^{pI}\nn\\
&=&P_I(-1)^{P_I(p-I)}e^I_j(-1)^{pI}=e_j^IP_I(-1)^{p(P_I+I)}~,\nn\eea
from this we see (\ref{temp1}) equals
\bea \textrm{(\ref{temp1})}=s^j\smalint_j~s^i\exp\Big\{e_j^I\frac{\partial}{\partial f_i^I}\Big\}\exp\Big\{e_j^I\frac{\partial}{\partial f_j^I}\Big\}=
s^is^j\smalint_j~\exp\Big\{e_j^I\big(\frac{\partial}{\partial f_i^I}+\frac{\partial}{\partial f_j^I}\big)\Big\}~.\label{concatenation}\eea
This is of course the discrete version of the familiar statement
\bea \int dx\int dy~\big(f(x)\delta(x-y)g(y)\big)=\int dx f(x)g(x)~,\nn\eea
and this analogy prompts us to define\bea \tilde m_{ij}=f_{iI}f^I_j(-1)^{pI}~,\nn\eea
and write (\ref{concatenation}) as
\bea \big[s^i\int\limits_i,\big[s^j\int\limits_j,\tilde m\big]\big]=s^is^j\int\limits_{i\cup j}~.\label{imagination}\eea
Thus $\tilde m_{ij}$ will serve as the '$\delta$ function'.

We also have the Stokes theorem
\bea [s^i\int\limits_i,\sum_k\,D^I_{~J}f_k^J\frac{\partial}{\partial f_k^I}]=\int\limits_i~D\exp\Big\{e_i^I\frac{\partial}{\partial f_i^I}\Big\}~,\nn\eea
which we write concisely as
\bea [s^i\int\limits_i,D]=s^i\int\limits_i d=0~,\label{Stokes}\eea
where on the left hand side  $D$ is an operator acting on the polynomials of $f$ while on the right hand side
 it is the differential of the Frobenius algebra.

Consider the polynomial in $s^i$, $s^{ij}$ which is the generating function of the graph cochains associated with the Frobenius algebra $\FR{a}$.
\bea \Gamma(\FR{a})=\exp\Big\{\sum\limits_i s^i\int\limits_i\Big\} \exp\Big\{\frac{1}{2}\sum_{i,j}f_i^I K_{IJ} f_j^{J}s^{ij}\Big\}\bigg|_{f=0}~.\label{master}\eea

We can deduce the properties of $\Gamma(\FR{a})$ by using the Stokes theorem.
Insert $D$ in between two exponentials,
\bea 0=P=\exp\Big\{\sum_i s^i\int\limits_i\Big\}\Big(\sum_iD^I_{~J}f_i^J\partial_{f_i^I}\Big) \exp\Big\{\frac12\sum_{i,j}f_i^I K_{IJ} f_j^{J}s^{ij}\Big\}\bigg|_{f=0}~.\label{master_P}\eea

Now instead of invoking Stokes theorem, we calculate $P$ by commuting $D$ to the rightmost position, doing so leaves us with only a commutator term
\bea \Big[D, \exp\Big\{\frac12\sum_{i,j}f_i^I K_{IJ} f_j^{J}s^{ij}\Big\}\Big]=\exp\big(\cdots\big)\frac12\Big\{f_i^LD^I_{~L}K_{IJ} f_j^{J}+
(-1)^{I}f_i^IK_{IL}D^L_{~J}f_j^{J}\Big\}s^{ij}~.\nn\eea
Naively we would conclude that the two terms in the braces cancel using the symmetry (\ref{sym_pairing}) and  (\ref{sym_K}).
 But one must remember in $(D^{-1})^I_{~J}$, the $I$ index must be exact while $J$ index is co-exact (or transverse to the exact part in general). Thus the first term in the curly brace is $-(-1)^Jm_{IJ}f_i^If_j^Js^{ij}$ but the $I$ index is exact while $J$ is
  co-exact; whereas the second term is $(-1)^Im_{IJ}f_i^If_j^Js^{ij}$, with $I$ co-exact and $J$ exact. To combine the two terms requires \emph{$p$ to be odd} since $\deg e^I+\deg e^J=p$ and we get
\bea
-\frac12\exp\big(\cdots\big)\sum_{i,j}(-1)^Jf_i^Im_{IJ}f_j^{J}s^{ij}=-\frac12\exp\big(\cdots\big)\sum_{i,j} \tilde m_{ij}s^{ij}~,\nn\eea
where the sum over $I,J$ indices are now over both exact and co-exact ones. Now we commute this term to the left most of (\ref{master_P}), again only picking up a commutator (using (\ref{imagination}))
\bea P=-\frac12\exp\Big\{\sum_k s^k\int\limits_k\Big\}\Big\{\sum_{i,j}s^is^js^{ij}\int\limits_{i\cup j}\Big\}\exp\Big\{\frac12\sum_{i,j}f_i^I K_{IJ} f_j^{J}s^{ij}\Big\}\bigg|_{f=0}~.\nn\eea
This term describes the splitting of a vertex into $i$ and $j$. In fact, we may write $P$ as
\bea &&P=-\frac12\Big(\sum_{pq}s^qs^{pq}S^p_q\Big)\exp\Big\{\sum_k s^k\int\limits_k\Big\}\exp\Big\{\frac12\sum_{i,j}f_i^I K_{IJ} f_j^{J}s^{ij}\Big\}\bigg|_{f=0}~,\nn\\
&&\hspace{2cm}S^p_q s^{ij}=\bigg\{\begin{array}{c}
                      s^{ij},~~i,j\neq p,q \\
                      s^{ip}+s^{iq},~~i\neq p,q
                    \end{array}~,\nn\eea
where we have defined a splitting operator $S$ acting on the polynomial of $s_{ij}$. And we have now found the differential operator for the graph cochains
\bea \delta_{Gph}=-\frac12\sum_{pq}s^qs^{pq}S^p_q~.\label{gph_codiff_poly}\eea
It is easy to see that this operator is nothing but the Fourier transform of the operator  (\ref{gph_diff_poly}).

To summarize, we have shown (\ref{master}) is a generating function for graph cocycles
\bea \delta_{Gph}\Gamma(\FR{a})=0~.\nn\eea

Next we show that the cohomology class of $\Gamma(\FR{a})$ is independent of the choice of $K$. The propagator is the combination of the inverse of $D$ and $m$:
$K_{IJ}=m_{IL}(D^{-1})^{L}_{~J}$. The inverse of $D$ is written with the help of Hodge decomposition. Under a change of the Hodge decomposition, $K$ changes by a $D$-exact term
\bea \gd K_{IJ}=J_{IL}D^L_{~J}+(-1)^{p-I}D^L_{~I}J_{LJ}~,~~~
J_{PQ}=(-1)^{PQ}J_{QP}~.\nn\eea
In particular
\bea\frac12\sum_{ij}f^I_i \delta K_{IJ} f^J_j=\frac12D\big((-1)^{p-1-I}J_{IJ}f_i^If_j^J\big)~.\nn\eea
Thus the corresponding change incurred in (\ref{master}) is
\bea \gd_K\Gamma(\FR{a})=\exp\Big\{\sum_i s^i\int\limits_i\Big\}\frac12D\big((-1)^{p-1-I}J_{IJ}f_i^If_j^J\big)\exp\Big\{\frac{1}{2}\sum_{i,j}f_i^I K_{IJ} f_j^{J}s^{ij}\Big\}\bigg|_{f=0}~,\nn\eea
Integrating by part, and we get
\bea\gd_K\Gamma(\FR{a})=-\exp\Big\{\sum_i s^i\int\limits_i\Big\}\frac12\big((-1)^{p-1-I}J_{IJ}f_i^If_j^J\big)\big(\frac12\sum_{ij}s^{ij}\tilde m_{ij}\big)\exp\Big\{\frac{1}{2}\sum_{i,j}f_i^I K_{IJ} f_j^{J}s^{ij}\Big\}\bigg|_{f=0}~,\nn\eea
and we manipulate $\tilde m$ in similar manner
\bea\gd_K\Gamma(\FR{a})&=&-\delta_{Gph}\big(\cdots\big)\nn\\
\big(\cdots\big)&=&\exp\Big\{\sum_i s^i\int\limits_i\Big\}\frac14\big((-1)^{p-1-I}J_{IJ}f_i^If_j^J\big)\exp\Big\{\frac{1}{2}\sum_{i,j}f_i^I K_{IJ} f_j^{J}s^{ij}\Big\}\bigg|_{f=0}~.\nn\eea

To summarize, we have shown that given an acyclic differential graded Frobenius algebra one can construct a class of graph cocycles. The explicit formula depends on the details of the propagator (Hodge decomposition),
but the change of propagator only causes the graph cycle to change by coboundaries and thus the class of cocycles is completely fixed by the data of the Frobenius algebra\footnote{Another proof using BV formalism was given in \cite{Hamilton-2007}, but the authors there used the Kontsevich theorem \ref{theorem-cocycle} as an input.} .
However, it is not clear to us how to remove the acyclicity condition.


\bigskip\bigskip

\begin{thebibliography}{6666}

\newcommand{\np}{{\em Nucl.\ Phys.\ }}
\newcommand{\pr}{{\em Phys.\ Rev.\ }}
\newcommand{\cmp}{{\em Commun.\ Math.\ Phys.\ }}
\newcommand{\pl}{{\em Phys.\ Lett.\ }}
%
%
\bibitem{Axelrod:1991vq}
  S.~Axelrod and I.~M.~Singer,
  ``Chern-Simons perturbation theory,''
  arXiv:hep-th/9110056;  ``Chern-Simons perturbation theory. 2,''
  J.\ Diff.\ Geom.\  {\bf 39} (1994) 173
  [arXiv:hep-th/9304087].
%
\bibitem{BN1}
D.~Bar-Natan,  ``On the Vassiliev Knot Invariants." Topology 34,
423-472, 1995.
\bibitem{Batalin:1981jr}
  I.~A.~Batalin and G.~A.~Vilkovisky,
  ``Gauge Algebra And Quantization,''
  Phys.\ Lett.\  B {\bf 102} (1981) 27.
%
\bibitem{Batalin:1984jr}
  I.~A.~Batalin and G.~A.~Vilkovisky,
  ``Quantization Of Gauge Theories With Linearly Dependent Generators,''
  Phys.\ Rev.\  D {\bf 28} (1983) 2567
  [Erratum-ibid.\  D {\bf 30} (1984) 508].
%
\bibitem{carmeli}
C.~Carmeli, L.~Caston and R.~Fioresi,
``Mathematical Foundation of Supersymmetry,''
  with an appendix with I. Dimitrov, {\it EMS Ser. Lect. Math.}, European Math. Soc., Zurich 2011.
%
\bibitem{Cattaneo:1999fm}
  A.~S.~Cattaneo and G.~Felder,
  ``A path integral approach to the Kontsevich quantization formula,''
  Commun.\ Math.\ Phys.\  {\bf 212} (2000) 591
  [arXiv:math/9902090].
%
\bibitem{cattaneo-GP}
A.~S.~Cattaneo, D.~Fiorenza, R.~Longoni,
``Graded Poisson Algebras,''
 {\it Encyclopedia of Mathematical Physics}, eds. J.-P. Franoise, G.L. Naber and Tsou S.T. ,
 vol. 2, p. 560-567 (Oxford: Elsevier, 2006).
%
\bibitem{Cattaneo:2008ph}
  A.~S.~Cattaneo and P.~Mn\"ev,
  ``Remarks on Chern-Simons invariants,''
  Commun.\ Math.\ Phys.\  {\bf 293} (2010) 803
  [arXiv:0811.2045 [math.QA]].
%
\bibitem{ConantVogtmann}
J.~Conant and K.~Vogtmann, ``On a Theorem of Kontsevich,'' Algebr.
Geom. Topol. 3 (2003) 1167-1224, arXiv:math/0208169.
%
\bibitem{susy-note-BB}
P.~Deligne, J.~W.~Morgan,
   ``Notes on supersymmetry (following Joseph Bernstein),"
 {\it Quantum fields and strings: a course for mathematicians}, Vol. 1, 2 (Princeton, NJ, 1996/1997),
 41Ð97, Amer. Math. Soc., Providence, RI, 1999.
%
\bibitem{Getzler:1994yd}
  E.~Getzler,
  ``Batalin-Vilkovisky algebras and two-dimensional topological field
  theories,''
  Commun.\ Math.\ Phys.\  {\bf 159} (1994) 265
  [arXiv:hep-th/9212043].

%
\bibitem{Hamilton-2005}
A.~Hamilton, ``A super-analogue of Kontsevich's theorem on graph
homology'', arXiv:math/0510390v1.
\bibitem{Hamilton-2007}
A.~Hamilton and A.~Lazarev, ''Graph cohomology classes in the Batalin-Vilkovisky formalism'', J.Geom.Phys.59:555-575,2009
%
\bibitem{HoschildSerre}
G.~Hochschild and J-P.~Serre, ''Cohomology of Lie Algebras'', The Annals of Mathematics, Second Series, Vol. 57, No. 3 (May, 1953), pp. 591-603
%
\bibitem{Kontsevich:symplectic}
M.~Kontsevich,  ``Formal (non)-commutative symplectic geometry,Ó The Gelfand Mathematical
Seminars, 1990 - 1992, Birkh\"auser (1993), 173 - 187.
%
\bibitem{Kontsevich:Feynamn}
  M.~Kontsevich, ``Feynman diagrams and low-dimensional topology,''
  First European Congress of Mathematics, 1992, Paris, Volume II,
  Progress in Mathematics {\bf 120}, Birkh\"auser 1994, 97 - 121.
%
\bibitem{Polyak:2004nn}
  M.~Polyak,
  ``Feynman diagrams for pedestrians and mathematicians,''
 {\it Graphs and patterns in mathematics and theoretical physics}, 15-42,
Proc. Sympos. Pure Math., {\bf 73}, Amer. Math. Soc., Providence, RI, 2005.
  [arXiv:math/0406251].
%
\bibitem{QiuZabzine:2009rf}
  J.~Qiu and M.~Zabzine,
  ``Odd Chern-Simons Theory, Lie Algebra Cohomology and Characteristic
  Classes,''~Commun.Math.Phys. {\bf 300}:789-833,2010,
  arXiv:0912.1243 [hep-th].
\bibitem{2010arXiv1006.1240Q}
J.~Qiu and M.~Zabzine, "Knot Invariants and New Weight Systems from General 3D TFTs", J.Geom.Phys. {\bf 62}:242-271, 2012, 
arXiv:1006.1240, hep-th. 
%
\bibitem{Roytenberg:2002nu}
  D.~Roytenberg,
  ``On the structure of graded symplectic supermanifolds and Courant
  algebroids,'' {\it Quantization, Poisson brackets and beyond (Manchester, 2001)}, 169-185,
   Contemp. Math., {\bf 315}, Amer. Math. Soc., Providence, RI, 2002
  [arXiv:math/0203110].
%
\bibitem{sawon1999}
J.~Sawon, ``Rozansky-Witten invariants of hyperk\"ahler
manifolds,''
 PhD thesis, Oxford 1999.
\bibitem{Sawon:2005qt}
  J.~Sawon,
  ``Perturbative expansion of Chern-Simons theory,''
  Geom. Topol. Monogr. {\bf 8} (2006) 145-166
  [arXiv:math/0504495].
%
\bibitem{Schwarz:1992nx}
  A.~S.~Schwarz,
  ``Geometry of Batalin-Vilkovisky quantization,''
  Commun.\ Math.\ Phys.\  {\bf 155} (1993) 249
  [arXiv:hep-th/9205088].
%
\bibitem{Schwarz:1999vn}
  A.~S.~Schwarz,
  ``Quantum observables, Lie algebra homology and TQFT,''
  Lett.\ Math.\ Phys.\  {\bf 49} (1999) 115
  [arXiv:hep-th/9904168].
%
\bibitem{Varadarajan}
V.~S.~Varadarajan,
 ``Supersymmetry for mathematicians: an introduction,''
 {\it Courant Lecture
Notes Series}, AMS, New York, 2004.
%
\end{thebibliography}
\end{document}